\documentclass[12pt,reqno]{amsart}

\usepackage[margin=1in]{geometry}
\usepackage{amsmath,amsfonts, amssymb, amsthm, mathtools, enumerate}
\usepackage{xcolor}
\usepackage{graphicx}
\usepackage{tikz}
\usepackage{tikz-cd}

\usepackage{hyperref}
\definecolor{lnkcol}{rgb}{0,0.25,0.75}
\hypersetup{colorlinks=true, linkcolor=lnkcol, citecolor=magenta, urlcolor=blue}

\newtheorem{theorem}{Theorem}[section]
\newtheorem{lemma}[theorem]{Lemma}
\newtheorem{proposition}[theorem]{Proposition}
\newtheorem{corollary}[theorem]{Corollary}
\theoremstyle{definition}
\newtheorem{definition}[theorem]{Definition}

\theoremstyle{remark}
\newtheorem{remark}[theorem]{Remark}

\numberwithin{figure}{section}
\numberwithin{equation}{theorem}

\title[Rigidity for Quasisymmetric Actions and $n$-Sphere Bowditch Boundaries]{Rigidity for Quasisymmetric Actions and $n$-Sphere Bowditch Boundaries} 
\author{Rana Sardar} 
\address{Department of Mathematical Sciences, Indian Institute of Science Education and Research Mohali, India} 
\email{ranasardar@iisermohali.ac.in} 
\subjclass[2020]{20F67, 20F65, 30L10, 57M60.} 
\keywords{Relatively hyperbolic groups, Bowditch boundaries, quasisymmetric maps, uniformization, Cannon conjecture}

\begin{document}
\begin{abstract}
    In this article, we establish a rigidity result for relatively hyperbolic groups whose Bowditch boundary is homeomorphic to a sphere, extending the main results of Bonk--Kleiner \cite{BK2002QM} and Kinneberg \cite{Kinneberg2015}. More precisely, let \((G,\mathcal P)\) be a relatively hyperbolic group pair whose Bowditch boundary \(Z\) is homeomorphic to an \(n\)-sphere and admits an Ahlfors \(Q\)-regular metric, where \(n\geq 2\). Under a suitable discrete path-length assumption on this metric, we prove that the action of \(G\) on \(Z\) is conjugate to an action by M\"obius transformations on the standard sphere \(\mathbb S^n\). As a consequence, \(G\) admits a properly discontinuous isometric action on \(\mathbb H^{n+1}\), and every peripheral subgroup is virtually isomorphic to \(\mathbb Z^n\).
\end{abstract}

\maketitle


\section{Introduction}
\label{sec:introduction}

A central theme in geometric group theory is the interplay between the large-scale geometry of a group and the geometry of its boundary. For Gromov-hyperbolic groups, the Gromov boundary provides a compact space that encodes important geometric information about the group. In particular, quasi-isometries between proper geodesic hyperbolic spaces induce quasisymmetric equivalences between their boundaries equipped with visual metrics. This connection makes the boundary a natural setting for studying geometric rigidity and uniformization problems. A fundamental question in this direction is whether the topology of the boundary, together with suitable metric and dynamical conditions, determines its quasisymmetric geometry. Of particular interest is the case where the boundary is homeomorphic to a sphere: under what conditions is it quasisymmetrically equivalent to the standard sphere? Furthermore, when does such a uniformization imply that the group admits a discrete, cocompact isometric action on real hyperbolic space?

\medskip

Let \(G\) be a non-elementary Gromov-hyperbolic group. Its boundary at infinity, denoted by \(\partial_\infty G\), is a compact metrizable space consisting of equivalence classes of geodesic rays. Visual metrics on \(\partial_\infty G\), arising from different visual parameters, belong to the same quasisymmetric gauge. Thus, the quasisymmetric geometry of the boundary is an invariant of the quasi-isometry type of \(G\). 
Relatively hyperbolic groups generalize hyperbolic groups by allowing distinguished peripheral subgroups. A relatively hyperbolic pair \((G,\mathcal P)\) admits a geometrically finite convergence action on a compact metrizable space, in which every limit point is either conical or bounded parabolic. The stabilizers of bounded parabolic points are precisely the conjugates of peripheral subgroups, acting properly discontinuously and cocompactly on the complement of the corresponding point. This compactum is called the Bowditch boundary and is denoted by \(\partial(G,\mathcal P)\). We refer to \cite{Gro87,Bowditch2012,Yaman2004,Hruska2010} for further details.

Kleinian groups provide fundamental examples of hyperbolic and relatively hyperbolic groups. A Kleinian group is a discrete subgroup \(\Gamma<\operatorname{PSL}(2,\mathbb C)\cong\operatorname{Isom}^{+}(\mathbb H^3)\), acting on \(\mathbb S^2\) by Möbius transformations. Its limit set \(\Lambda(\Gamma)\) consists of accumulation points of an orbit in \(\mathbb H^3\). If \(\Gamma\) acts cocompactly on \(\mathbb H^3\), it is a hyperbolic group with boundary equivariantly homeomorphic to \(\mathbb S^2\). Such a group is a uniform lattice, and its quotient is a compact hyperbolic \(3\)-orbifold. More generally, geometrically finite Kleinian groups with parabolic elements are relatively hyperbolic with respect to their maximal parabolic subgroups, and their Bowditch boundaries are equivariantly homeomorphic to their limit sets. In particular, finite-volume hyperbolic \(3\)-manifolds with cusps provide examples with Bowditch boundary homeomorphic to \(\mathbb S^2\).

The relationship between boundary topology and geometric realizations in hyperbolic space lies at the heart of Cannon-type conjectures. In dimension one, the situation is well understood: every hyperbolic group with boundary homeomorphic to \(\mathbb S^1\) is virtually a closed surface group and acts geometrically on \(\mathbb H^2\). In dimension two, Cannon's conjecture asserts that every hyperbolic group with boundary homeomorphic to \(\mathbb S^2\) admits a proper, cocompact, isometric action on \(\mathbb H^3\), possibly with finite kernel \cite[Conjecture~9.1]{KB2002}, \cite[Problem~53]{Kapovich2005}. 
Replacing the Gromov boundary with the Bowditch boundary leads to the relative Cannon conjecture \cite[Problem~57]{Kapovich2005}. This generalization is motivated by finite-volume hyperbolic \(3\)-manifolds with cusps, whose fundamental groups are relatively hyperbolic with respect to cusp subgroups, which are virtually \(\mathbb Z^2\), and have Bowditch boundary homeomorphic to \(\mathbb S^2\).

A fruitful approach to Cannon-type conjectures is through the metric geometry of group boundaries, particularly quasisymmetric geometry. Bonk and Kleiner proved that a compact Ahlfors \(Q\)-regular Loewner space homeomorphic to \(\mathbb S^2\) has \(Q=2\) and is quasisymmetrically equivalent to the standard sphere \cite[Theorem~1.2]{BK2002Qs}. Using this uniformization theorem, they obtained a conditional solution to Cannon's conjecture: if a hyperbolic group has boundary homeomorphic to \(\mathbb S^2\) and its Ahlfors-regular conformal dimension is attained, then it acts geometrically on \(\mathbb H^3\) \cite{BK2005}. 

In the relatively hyperbolic setting, however, additional difficulties arise because Bowditch boundaries need not be doubling and may not admit Ahlfors-regular metrics. Several works have addressed rigidity and uniformization in this setting, including those of Markovic, Haïssinsky, Einstein, Groves, Tshishiku, Walsh, Hruska, and many others \cite{Mar2013,Hai2015,GE2020,TW2020,HW2023,HW2024,GHMOSW2026}. In particular, Mackay and Sisto established a connection between the doubling property of Bowditch boundaries and the algebraic structure of peripheral subgroups, under suitable finite generation assumptions \cite{MS2020}. More recently, Pal and Sardar proved a relative counterpart of the Bonk--Kleiner theorem: if a Bowditch boundary homeomorphic to \(\mathbb S^2\) admits an Ahlfors \(Q\)-regular metric realizing its Ahlfors-regular conformal dimension, then the relatively hyperbolic group admits a discrete isometric action on \(\mathbb H^3\), with peripheral subgroups virtually \(\mathbb Z^2\) \cite{PS2026}.

\medskip
The higher-dimensional situation is markedly different. For \(n\geq3\), a direct analogue of Cannon's conjecture is false in general: there exist hyperbolic groups with Gromov boundary homeomorphic to \(\mathbb S^n\) that do not admit geometric actions on real hyperbolic space \(\mathbb H^{n+1}\). The Gromov--Thurston examples of closed negatively curved manifolds that do not admit hyperbolic metrics provide such examples \cite{GromovThurston1987}. Another family is furnished by closed complex hyperbolic \(2\)-manifolds. Their fundamental groups are hyperbolic and have Gromov boundary homeomorphic to \(\mathbb S^3\). Moreover, this boundary admits an Ahlfors \(4\)-regular metric realizing its Ahlfors-regular conformal dimension \cite{Mackay-Tyson}, yet the group cannot act discretely and cocompactly on real hyperbolic \(4\)-space. Nevertheless, spherical boundary topology continues to impose strong restrictions on hyperbolic groups. Bartels, Lück, and Weinberger proved that every torsion-free hyperbolic group with boundary homeomorphic to \(\mathbb S^n\), for \(n\geq5\), is the fundamental group of a closed aspherical \(n\)-manifold \cite{BartelsLuckWeinberger2010}. More recently, Bregman and Incerti-Medici obtained related manifold-realization results for certain cubulated hyperbolic groups with spherical boundaries, under additional hypotheses \cite{BregmanIncertiMedici2026}. These examples and results demonstrate that boundary topology alone does not determine a real hyperbolic uniformization in higher dimensions. They instead point to the importance of additional metric and dynamical conditions.


One approach to higher-dimensional rigidity and uniformization is to study the dynamics of group actions on metric spaces. In this direction, Bonk and Kleiner established a uniformization theorem for uniformly quasi-M\"obius actions on compact Ahlfors \(n\)-regular metric spaces of topological dimension \(n\). They showed that if the induced action on the space of distinct triples is cocompact, then the underlying space is quasisymmetrically equivalent to the standard sphere, and the group action is conjugate to an action by M\"obius transformations \cite{BK2002QM}. Their theorem reveals a fundamental connection between the dynamics of group actions and the metric uniformization of spheres. The cocompactness condition on triples is particularly natural for hyperbolic groups, whose boundary actions are uniform convergence actions.
Kinneberg subsequently developed this approach in \cite[Theorem~1.5]{Kinneberg2015}. In particular, he established a de-snowflaking theorem for Ahlfors \(Q\)-regular metric \(n\)-spheres \((Z,d)\) equipped with strongly quasi-M\"obius group actions that are cocompact on triples and satisfy an additional discrete path lower-bound condition. Under these hypotheses, he constructed a metric \(d_{\mathrm{new}}\) on \(Z\) that is bi-Lipschitz equivalent to the standard spherical metric. Consequently, \(Z\) is quasisymmetrically equivalent to \(\mathbb S^n\). Furthermore, when \(n\geq 2\), the group action is conjugate to an action by M\"obius transformations.

\medskip
\paragraph{\textbf{Main Results}}

The preceding developments provide the motivation for the present work. We investigate metric uniformization mechanisms in the context of relatively hyperbolic groups, where the action on the Bowditch boundary is naturally studied through its induced action on the space of distinct pairs. Our goal is to extend the metric uniformization techniques of Bonk--Kleiner and Kinneberg to this relative hyperbolic setting.

Our first main result establishes a uniformization theorem for pair-cocompact actions. It shows that, under suitable metric and dynamical hypotheses, Ahlfors regularity and the topological dimension of the underlying space determine its quasisymmetric equivalence class.

\begin{theorem}\label{thm_main_1}
    Let \(n\in\mathbb N\), and let \(Z\) be a compact, Ahlfors \(n\)-regular metric space of topological dimension \(n\). Suppose that \(G\curvearrowright Z\) is a uniformly quasisymmetric action and that the induced action $G\curvearrowright \operatorname{Pair}(Z)$ is cocompact. Then \(Z\) is quasisymmetrically equivalent to the standard sphere \(\mathbb S^n\).
\end{theorem}

Our second main result is a de-snowflaking theorem. It provides a method for modifying the metric on a compact metric sphere so that its Ahlfors regular dimension agrees with its topological dimension. The key ingredient is a discrete length condition, which allows us to construct a metric in the same quasisymmetric gauge as the original one. 
Recall that, for \(\tau>0\), a discrete \(\tau\)-path joining \(x\) to \(y\) is a finite sequence \(x=z_0,z_1,\ldots,z_\ell=y\) such that \(d(z_{i-1},z_i)\leq\tau\) for every \(i\). The integer \(\ell\) is called the combinatorial length of the path.

\begin{theorem}\label{thm_main_2}
    Let \(n\in\mathbb N\), \(0<\epsilon\leq 1\), and let \((Z,d)\) be a metric space homeomorphic to \(\mathbb S^n\) and Ahlfors \(Q\)-regular, where \(Q=n/\epsilon\). Suppose that \(G\curvearrowright Z\) is a uniformly power quasisymmetric action that is cocompact on pairs. Assume, moreover, that there exists a constant \(c>0\) such that every discrete \(\tau\)-path joining \(x\) to \(y\) satisfies
    \[
    \ell\geq c\left(\frac{d(x,y)}{\tau}\right)^{1/\epsilon}.
    \]
    Then there exist a metric \(d_{\mathrm{new}}\) on \(Z\) and a constant \(C_*\geq1\) such that
    \[
    C_*^{-1}d(x,y)^{1/\epsilon} \leq d_{\mathrm{new}}(x,y) \leq C_*d(x,y)^{1/\epsilon}
    \]
    for all \(x,y\in Z\).
\end{theorem}

Combining Theorems~\ref{thm_main_1} and \ref{thm_main_2}, we obtain the following uniformization result.

\begin{corollary}\label{cor:thm-2}
    Under the hypotheses of Theorem~\ref{thm_main_2}, there exists a metric \(d_{\mathrm{new}}\) on \(Z\) such that \((Z,d_{\mathrm{new}})\) is Ahlfors \(n\)-regular and quasisymmetrically equivalent to the standard sphere \(\mathbb S^n\). Inparticular, $(Z,d)$ quasisymmetrically equivalent to the standard sphere \(\mathbb S^n\). 
\end{corollary}

\begin{proof}
    
    By Theorem~\ref{thm_main_2}, there exists a metric \(d_{\mathrm{new}}\) on \(Z\) such that
    \[
    d_{\mathrm{new}}(x,y)\asymp d(x,y)^{1/\epsilon}.
    \]
    Consequently, the identity map $\operatorname{id}:(Z,d)\longrightarrow (Z,d_{\mathrm{new}})$ is quasisymmetric. Moreover, since \(d_{\mathrm{new}}\) is comparable to a positive power of \(d\), the two metrics induce the same topology on \(Z\).

    Let \(\mu\) be an Ahlfors \(Q\)-regular measure on \((Z,d)\). The metric comparison above implies that, for every \(x\in Z\) and sufficiently small \(r>0\),
    \[
    \mu(B_{d_{\mathrm{new}}}(x,r)) \asymp \mu(B_d(x,r^\epsilon)) \asymp (r^\epsilon)^Q =r^n.
    \]
    Thus \((Z,d_{\mathrm{new}})\) is Ahlfors \(n\)-regular. Since the two metrics induce the same topology, we also have $\dim_{\mathrm{top}}(Z,d_{\mathrm{new}})=n$.

    Furthermore, conjugating the original action by the identity map gives a uniformly quasisymmetric action $G\curvearrowright (Z,d_{\mathrm{new}})$. Since the two metrics induce the same topology, the corresponding spaces of distinct pairs are identical as topological spaces. In particular, the cocompactness of the action on pairs is preserved.

    We can therefore apply Theorem~\ref{thm_main_1} to \((Z,d_{\mathrm{new}})\). It follows that \((Z,d_{\mathrm{new}})\) is quasisymmetrically equivalent to \(\mathbb S^n\), and hence \((Z,d)\) is quasisymmetrically equivalent to \(\mathbb S^n\), as required.
\end{proof}

We finally apply these metric uniformization results to Bowditch boundaries of relatively hyperbolic groups, using the fact that every relatively hyperbolic group acts on its Bowditch boundary by uniformly power-quasisymmetric homeomorphisms. Under additional metric assumptions on the boundary, we obtain a higher-dimensional analogue of Cannon-type uniformization.

\begin{theorem}\label{thm:generalized-rel-Cannon-conjec}
    Let \((G,\mathcal P)\) be a relatively hyperbolic group pair whose Bowditch boundary \(Z=\partial(G,\mathcal P)\) is homeomorphic to \(\mathbb S^n\), where \(n\geq2\). Suppose that \(Z\) admits an Ahlfors \(Q\)-regular metric \(d\), where \(Q=n/\epsilon\) for some \(0<\epsilon\leq1\), and that there exists a constant \(c>0\) such that every discrete \(\tau\)-path joining \(x\) to \(y\) has combinatorial length satisfying
    \[
    \ell\geq c\left(\frac{d(x,y)}{\tau}\right)^{1/\epsilon}.
    \]
    Then the action of $G$ on $Z$ is conjugate to an action by M\"obius transformations on $\mathbb S^n$. In particular, $G$ admits a properly discontinuous isometric action on $\mathbb H^{n+1}$, and every peripheral subgroup is virtually isomorphic to $\mathbb Z^n$.
\end{theorem} 

\medskip
\paragraph{\textbf{Organization}} 
The paper is organized as follows. Section 2 collects the necessary preliminaries on relatively hyperbolic groups, metric properties, dimensions, and quasi-M\"obius and quasisymmetric maps. Sections 3–5 develop the main tools concerning quasisymmetric dynamics, weak tangents, stable values, bounded multiplicity, and regular maps. Section 6 is devoted to the proof of Theorem~\ref{thm_main_1}. In Section 7, we introduce a de-snowflaking procedure and prove Theorem~\ref{thm_main_2}. Finally, Section 8 contains the proof of Theorem~\ref{thm:generalized-rel-Cannon-conjec}.

\section{Preliminaries}

\subsection{Relatively Hyperbolic Groups and Boundaries}
\leavevmode 
\smallskip

Several notions of relative hyperbolicity appear in the literature. We refer to Hruska’s article \cite{Hruska2010} for a discussion of the equivalence of these definitions.
In this article,  we adopt the dynamical characterization of relative hyperbolicity due to Gerasimov \cite[Definition $RH_{32}$]{Gera2012} as our definition. 

\begin{definition}[Relatively hyperbolic group] \label{defn:rel_hyp_gp}
    Let $G$ be a finitely generated discrete group, and let $\mathcal{H}=\{H_1,\ldots,H_m\}$ be a finite collection of subgroups of $G$. We say that the pair $(G,\mathcal{H})$ is \emph{relatively hyperbolic} if there exists a compact Hausdorff space $Z$ together with an action \(G\curvearrowright Z\) by homeomorphisms satisfying the following conditions:
    \begin{enumerate}[$(i)$]
        \item The action is a \emph{$3$-proper} (or \emph{convergence action}); that is, the induced action on the space of distinct ordered triples $\operatorname{Tri}(Z)$ is proper, where
        \[
        \operatorname{Tri}(Z)=\{(x,y,z)\in Z^3 : x,y,z \text{ are pairwise distinct}\}.
        \]

        \item The action is \emph{$2$-cocompact} (or \emph{expansive}); namely, \(\operatorname{Pair}(Z)/G\) is compact, where
    \[
    \operatorname{Pair}(Z)=\{(x,y)\in Z\times Z : x\neq y\}.
    \]

        \item The collection $\mathcal{H}$ consists of representatives of the conjugacy classes of stabilizers of the parabolic points of the action.
    \end{enumerate}

    Here, a point $p\in Z$ is called \emph{parabolic} if $\operatorname{Stab}_G(p)$ is infinite and acts cocompactly on $Z\setminus\{p\}$. The subgroup $\operatorname{Stab}_G(p)$ is called the \emph{parabolic subgroup} associated to $p$.
\end{definition} 

The compactum $Z$ in Definition~\ref{defn:rel_hyp_gp} is called  the Bowditch boundary of $(G,\mathcal H)$. We denote it by $\partial(G,\mathcal H)$. In Definition~\ref{defn:rel_hyp_gp}, the induced action $G\curvearrowright\operatorname{Tri}(Z)$ is cocompact if and only if the group $G$ is hyperbolic. 
Moreover, since the natural projection $\operatorname{Tri}(Z)\to \operatorname{Pair}(Z)$ defined by $(x,y,z)\mapsto (x,y)$, is continuous and surjective, every $3$-cocompact action is also $2$-cocompact.

\subsection{Metric Properties and Dimensions}
\label{subsec:metric_properties}

\begin{definition}[Uniformly perfect metric space] 
    A metric space $(Z,d)$ is \emph{uniformly perfect} if there exists $C>1$ such that, for every $a\in Z$ and every $0<R\leq\operatorname{diam}(Z)$, we have $B(a,R)\setminus B\left(a,\frac{R}{C}\right)\neq\varnothing$. 
\end{definition} 

\begin{definition}[Doubling metric space] 
    A metric space $(Z,d)$ is said to be \emph{$N$-doubling} if every ball of radius $R$ can be covered by at most $N$ balls of radius $R/2$. 
\end{definition} 

The boundary of a hyperbolic group is doubling. In contrast, the Bowditch boundary of a relatively hyperbolic group need not be doubling. The following result characterizes doubling of the Bowditch boundary in terms of the peripheral subgroups.

\begin{proposition}[{\cite[Proposition~4.5]{MS2020}}] 
    Let $(G,\mathcal H)$ be relatively hyperbolic. Then its Bowditch boundary is doubling if and only if every peripheral subgroup in $\mathcal H$ is virtually nilpotent. 
\end{proposition}

\begin{definition}[Ahlfors regular metric measure space] \label{defn:Ahlfors_reg} 
    Let $(Z,d,\mu)$ be a metric measure space. For $Q>0$, we say that $(Z,d,\mu)$ is \emph{Ahlfors $Q$-regular} if there exists a constant $C\geq1$ such that for every $z\in Z$ and every $0<r\leq\operatorname{diam}(Z)$, we have
    \[ 
    C^{-1}r^Q \leq \mu(B(z,r)) \leq Cr^Q.
    \]     
\end{definition} 

By a result of Coornaert \cite{Coornaert1993}, the Gromov boundary of a hyperbolic group is Ahlfors regular with respect to a visual metric. In contrast, for a relatively hyperbolic group, a visual metric on the Bowditch boundary need not be Ahlfors regular; \cite[Proposition~1.6]{BM1996}, \cite[Example~1, p.~189]{GP2001}. 


\begin{definition}[Ahlfors regular conformal dimension] \label{defn:conf_dim}
    Let $(Z,d)$ be a metric space. The \emph{Ahlfors regular conformal dimension} of $Z$ is defined by
    \[
    \operatorname{Cdim}_{\mathrm{AR}}(Z) = \inf\{\dim_H(X): X \text{ is Ahlfors regular and quasisymmetric to } Z\}.
    \]
\end{definition}

 
\subsection{Quasi-M\"{o}bius and Quasisymmetric Maps}

\begin{definition}
    Let $(Z,d_Z)$ be a metric space. The \emph{cross-ratio} \([z_1,z_2,z_3,z_4]\) of four distinct points $z_1,z_2,z_3,z_4$ in $Z$ is defined by
    \[
    [z_1,z_2,z_3,z_4]
    :=
    \frac{d_Z(z_1,z_3)\, d_Z(z_2,z_4)}
     {d_Z(z_1,z_4)\, d_Z(z_2,z_3)}.
    \]
    
    Let $\eta:[0,\infty)\to[0,\infty)$ be a homeomorphism, and let \(f:(X,d_X)\to(Y,d_Y)\) be an injective map between metric spaces. 
    The map $f$ is called \emph{$\eta$-quasi-M\"{o}bius} if for every four-tuple $(x_1,x_2,x_3,x_4)$ of distinct points in $X$,
    \[
    [f(x_1),f(x_2),f(x_3),f(x_4)]
    \leq
    \eta\,\bigl([x_1,x_2,x_3,x_4]\bigr).
    \]
    The map $f$ is called \emph{$\eta$-quasisymmetric} if for every triple $(x_1,x_2,x_3)$ of distinct points in $X$,
    \[
    \frac{d_Y\bigl(f(x_1),f(x_2)\bigr)}
     {d_Y\bigl(f(x_1),f(x_3)\bigr)}
    \leq \eta\!\left( \frac{d_X(x_1,x_2)} {d_X(x_1,x_3)} \right). 
    \]

    A quasisymmetric map is called \emph{power quasisymmetric} if its distortion function has the form
    \[
    \eta(t) =
    \begin{cases}
        C t^{1/\alpha} & \text{for } 0 < t < 1, \\[4pt]
        C t^{\alpha} & \text{for } 1 \le t.
    \end{cases}
    \]
    for some $C \geq 1$ and $\alpha>0$. This is also called a power-type quasisymmetric map. 
\end{definition}

Note that Quasi-M\"{o}bius and quasisymmetric maps are homeomorphisms onto their images. The composition of two quasi-M\"{o}bius maps with distortion functions $\eta_1$ and $\eta_2$ is $\eta_2\circ\eta_1$-quasi-M\"{o}bius. Similarly, compositions of quasisymmetric maps are quasisymmetric, and the same holds for power quasisymmetric maps. Moreover, every $\eta$-quasisymmetric map is $\widetilde{\eta}$-quasi-M\"{o}bius, where $\widetilde{\eta}$ depends only on $\eta$. Conversely, every quasi-M\"{o}bius map between bounded metric spaces is quasisymmetric. In particular, for the Gromov or Bowditch boundaries, the classes of power quasisymmetric, quasisymmetric, and quasi-M\"obius homeomorphisms coincide \cite{Pau1996,Bonk2000,MS2024,PS2024,Sardar2026}. 
A power quasisymmetric map with distortion function
\[
\eta(t)=
\begin{cases}
Ct^{1/\alpha}, & 0<t<1,\\
Ct^\alpha, & t\geq 1,
\end{cases}
\qquad C\geq1,\quad \alpha>0,
\] 
is also a quasisymmetric map with distortion function
\[
\eta'(t)=C\max\{t^{1/\alpha},t^\alpha\}.
\]

\section{Cocompactness on Pairs and Quasisymmetric Dynamics}
\label{sec:cocompact-actions}

In this section, we study consequences of the cocompactness of the induced action on the space of distinct pairs. We first give a quantitative characterization of this cocompactness and then establish several technical lemmas that will be used in subsequent sections. Throughout, $B(p,r)$ denotes the open ball of radius $r$ centered at $p$.

Let $(Z,d_Z)$ be a compact metric space. The space of distinct ordered pairs in $Z$ is
\[
\operatorname{Pair}(Z)
=\{(z_1,z_2)\in Z\times Z:z_1\ne z_2\}.
\]
An action $G\curvearrowright Z$ by homeomorphisms induces an action on $\operatorname{Pair}(Z)$ by $g\cdot(z_1,z_2)=(gz_1,gz_2)$. We say that this induced action is \emph{cocompact} if there is a nonempty compact set $K\subset\operatorname{Pair}(Z)$ such that
\[
\operatorname{Pair}(Z)=GK=\bigcup_{g\in G}gK.
\]
The following lemma translates this topological condition into a uniform separation property.

\begin{lemma}[Uniform separation characterization]
\label{lem:cocompact_uniformly_separated}
    Assume that $Z$ contains at least two points. The action $G\curvearrowright\operatorname{Pair}(Z)$ is cocompact if and only if there exists $\delta>0$ such that, for every $(z_1,z_2)\in\operatorname{Pair}(Z)$, there is an element $g\in G$ satisfying $d_Z(gz_1,gz_2)\geq\delta$.
\end{lemma}

\begin{proof}
    Suppose first that the action on $\operatorname{Pair}(Z)$ is cocompact. Then there is a nonempty compact set $K\subset\operatorname{Pair}(Z)$ such that $\operatorname{Pair}(Z)=GK$. Since the distance function is continuous and positive on $K$, compactness gives
    \[
    \delta:=\min\{d_Z(x,y):(x,y)\in K\}>0.
    \]
    Given a pair $(z_1,z_2)$ of distinct points, choose $h\in G$ such that $(z_1,z_2)\in hK$. Then $(h^{-1}z_1,h^{-1}z_2)\in K$, so
    \[
    d_Z(h^{-1}z_1,h^{-1}z_2)\geq\delta.
    \]
    Taking $g=h^{-1}$ proves the required separation property.

    Conversely, suppose that such a $\delta>0$ exists, and set
    \[
    K_\delta:=\{(z_1,z_2)\in Z\times Z:d_Z(z_1,z_2)\geq\delta\}.
    \]
    The set $K_\delta$ is closed in the compact space $Z\times Z$, hence compact. Since $\delta>0$, it is contained in $\operatorname{Pair}(Z)$. By assumption, every pair in $\operatorname{Pair}(Z)$ can be mapped into $K_\delta$ by an element of $G$. Equivalently, every pair belongs to some translate of $K_\delta$. Thus
\[
\operatorname{Pair}(Z)=GK_\delta,
\]
which proves cocompactness.
\end{proof}

The next lemma is a rescaling principle for quasisymmetric homeomorphisms of a uniformly perfect compact metric space. We use the standard fact that a family of uniformly quasisymmetric maps is also uniformly quasi-M\"obius; the associated quasi-M\"obius control function depends only on the quasisymmetric control function. 

\begin{lemma} \label{lem:tech_2}
    Let $(Z,d_Z)$ be a uniformly perfect compact metric space. Let $\delta>0$, and let $\eta:(0,\infty)\to(0,\infty)$ be a quasisymmetry control function. For each $k\in\mathbb N$, suppose we are given
    \begin{itemize}
        \item a ball $B_k=B(p_k,R_k)\subseteq Z$;
        
        \item distinct points $x_k^1,x_k^2\in B(p_k,\lambda_kR_k)$ such that $d_Z(x_k^1,x_k^2)>\delta_kR_k$, where $\lambda_k,\delta_k>0$ and $\lambda_k\to0$, $\delta_k\to0$;
        
        \item an $\eta$-quasisymmetric homeomorphism $g_k:Z\to Z$ such that, writing $y_k^i:=g_k(x_k^i)$, we have $d_Z(y_k^1,y_k^2)>\delta$.
    \end{itemize}
    Then the following statements hold:
    \begin{enumerate}
        \item If $D_k\subseteq B_k$ is an $\varepsilon_kR_k$-dense subset of $B_k$, where $\varepsilon_k>0$ and $\varepsilon_k/\delta_k\to0$, then
        \[
        \operatorname{dist}_H\bigl(g_k(D_k),Z\bigr)\longrightarrow0.
        \]

        \item The image of $B_k$ becomes asymptotically dense in $Z$, in the sense that
        \[
        \operatorname{diam}\bigl(Z\setminus g_k(B_k)\bigr)\longrightarrow0.
        \]
    \end{enumerate}
\end{lemma}

\begin{proof}
    Let $\eta_{\mathrm{qm}}$ be a quasi-M\"obius control function for the maps $g_k$, chosen independently of $k$. Such a function exists because the maps are uniformly $\eta$-quasisymmetric.

    \medskip
    \noindent\textbf{(1)} Since $g_k(D_k)\subseteq Z$, it is enough to show that every point of $Z$ is at a distance tending uniformly to zero from $g_k(D_k)$. Fix $z\in Z$ and put $z_k:=g_k^{-1}(z)$.

    \smallskip
    \noindent\emph{Case 1.} Suppose that $z_k\in B_k$. Choose $y_k\in D_k$ such that
    \[
    d_Z(y_k,z_k)\leq\varepsilon_kR_k.
    \]
    Since $d_Z(x_k^1,x_k^2)>\delta_kR_k$, at least one point $a_k\in\{x_k^1,x_k^2\}$ satisfies
    \[
    d_Z(y_k,a_k)\geq\frac{\delta_kR_k}{2}.
    \]
    For all sufficiently large $k$, the condition $\varepsilon_k/\delta_k\to0$ ensures that the points in the following quasisymmetry comparison are distinct (the conclusion is immediate if $y_k=z_k$). Hence
    \begin{align*}
        d_Z\bigl(g_k(y_k),z\bigr)
        &\leq d_Z\bigl(g_k(y_k),g_k(a_k)\bigr)
       \eta\!\left(\frac{d_Z(y_k,z_k)}{d_Z(y_k,a_k)}\right)\\
        &\leq \operatorname{diam}(Z)
       \eta\!\left(\frac{2\varepsilon_k}{\delta_k}\right).
    \end{align*}
    This bound tends to zero.

    \smallskip
    \noindent\emph{Case 2.} Suppose that $z_k\notin B_k$. Then $d_Z(z_k,p_k)\geq R_k$. Uniform perfectness gives a constant $C>0$, independent of $k$, and a point $w_k\in B_k$ such that
    \[
    d_Z(w_k,p_k)\geq CR_k.
    \]
    Since $D_k$ is $\varepsilon_kR_k$-dense in $B_k$, choose $y_k\in D_k$ with $d_Z(y_k,w_k)\leq\varepsilon_kR_k$. Thus, $d_Z(y_k,p_k)\geq (C-\varepsilon_k)R_k$. 
    The assumptions $\delta_k \to 0$ and $\varepsilon_k/\delta_k\to0$ imply $\varepsilon_k\to0$. Thus, for all sufficiently large $k$ we have
    \[
    d_Z(y_k,p_k)\geq C_0R_k,
    \qquad C_0:=C/2>0.
    \]
    We may also assume that $\lambda_k<\min\{1/2,C_0/2\}$. Using the quasi-M\"obius inequality for the quadruple $(z_k,y_k,x_k^1,x_k^2)$, we obtain
    \begin{align*}
        \frac{d_Z\bigl(g_k(z_k),g_k(y_k)\bigr)d_Z(y_k^1,y_k^2)}{d_Z\bigl(g_k(z_k),y_k^1\bigr)d_Z\bigl(g_k(y_k),y_k^2\bigr)}
        & \leq \eta_{\mathrm{qm}}\!\left( \frac{d_Z(z_k,y_k)d_Z(x_k^1,x_k^2)}{d_Z(z_k,x_k^1)d_Z(y_k,x_k^2)}\right)\\
        & \leq \eta_{\mathrm{qm}}\!\left( \frac{4d_Z(z_k,p_k)\lambda_k}{(d_Z(z_k,p_k)-\lambda_kR_k)(C_0-\lambda_k)} \right)\\
        & \leq \eta_{\mathrm{qm}}\!\left(\frac{16\lambda_k}{C_0}\right).
    \end{align*}
    In the last inequality we used $d_Z(z_k,p_k)\geq R_k$, $\lambda_k\leq1/2$, and $C_0-\lambda_k\geq C_0/2$. Since $d_Z(y_k^1,y_k^2)>\delta$ and the other two image distances in the denominator are at most $\operatorname{diam}(Z)$, it follows that
    \[
    d_Z\bigl(g_k(z_k),g_k(y_k)\bigr)
    \leq\frac{\operatorname{diam}(Z)^2}{\delta} \eta_{\mathrm{qm}}\!\left(\frac{16\lambda_k}{c_0}\right) \longrightarrow0.
    \]
    The estimates in the two cases are uniform in $z$. This proves (1).

    \medskip
    \noindent\textbf{(2)} Fix $k$ and let $y_k',z_k'\in Z\setminus g_k(B_k)$. If the two points coincide, there is nothing to prove. Otherwise, set
    \[
    y_k:=g_k^{-1}(y_k'),\qquad z_k:=g_k^{-1}(z_k').
    \]
    Then $y_k,z_k\notin B_k$. After interchanging their roles if necessary, assume
    \[
    d_Z(y_k,p_k)\leq d_Z(z_k,p_k).
    \]
    For all sufficiently large $k$, we have $\lambda_k<1/2$. The quasi-M\"obius inequality gives
    \begin{align*}
        \frac{d_Z(z_k',y_k')d_Z(y_k^1,y_k^2)}{d_Z(z_k',y_k^1)d_Z(y_k',y_k^2)}
        & \leq \eta_{\mathrm{qm}}\!\left(\frac{d_Z(z_k,y_k)d_Z(x_k^1,x_k^2)}{d_Z(z_k,x_k^1)d_Z(y_k,x_k^2)}\right)\\
        & \leq\eta_{\mathrm{qm}}\!\left(\frac{4d_Z(z_k,p_k)\lambda_k}{(d_Z(z_k,p_k)-\lambda_kR_k)(1-\lambda_k)}\right)\\
        & \leq\eta_{\mathrm{qm}}(16\lambda_k).
    \end{align*}
    Since $d_Z(y_k^1,y_k^2)>\delta$ and the other two image distances in the denominator are at most $\operatorname{diam}(Z)$, we conclude that
    \[
    d_Z(y_k',z_k')
    \leq\frac{\operatorname{diam}(Z)^2}{\delta} \eta_{\mathrm{qm}}(16\lambda_k)\longrightarrow0.
    \]
    The bound is independent of the choice of $y_k',z_k'$. Therefore
    \[
    \operatorname{diam}\bigl(Z\setminus g_k(B_k)\bigr)\longrightarrow0,
    \]
    as required.
\end{proof}


We now specialize to uniformly power-quasisymmetric actions on metric spheres. Let $(Z,d_Z)$ be a metric space homeomorphic to the standard sphere $\mathbb S^n$, with $n\geq1$. Suppose that $G\curvearrowright Z$ is an action by uniformly power-quasisymmetric homeomorphisms with common distortion function
\[
\eta(t)=C\max\{t^\alpha,t^{1/\alpha}\},
\qquad C\geq1,\quad \alpha>0.
\]
Assume also that the induced action on $\operatorname{Pair}(Z)$ is cocompact. By the standard implication from quasisymmetry to quasi-M\"obius distortion, there is a quasi-M\"obius control function $\eta_{\mathrm{qm}}$, depending only on $C$ and $\alpha$, that works for every element of $G$. 
Throughout the remainder of this section, for a fixed $g\in G$ we write $x'=g(x)$ for $x\in Z$.

\begin{lemma}
\label{lem:power_quasisymmetric_estimates}
    There exists constant $C_1\geq1$, depending only on $C,\alpha$ and the constant $\delta$ considered in Lemma~\ref{lem:cocompact_uniformly_separated} such that the following holds. Given $p\in Z$, $0<r<\operatorname{diam}(Z)$, and $\lambda\geq2$, there is an element $g\in G$ depending on $p$, $r$ and $\operatorname{diam}(Z)$ such that we have:
    \begin{enumerate}
        \item For every pair of distinct points $x,y\in Z$,
        \[
        d_Z(x',y')\leq C_1\max\left\{\left(\frac{d_Z(x,y)}{r}\right)^\alpha,\left(\frac{d_Z(x,y)}{r}\right)^{1/\alpha}\right\}.
        \]
        
        \item For every pair of distinct points $x,y\in Z$,
        \[
        \frac{r}{d_Z(x',y')}\leq C_1\max\left\{\left(\frac{\operatorname{diam}(Z)}{d_Z(x,y)}\right)^\alpha,\left(\frac{\operatorname{diam}(Z)}{d_Z(x,y)}\right)^{1/\alpha}\right\}.
        \]
        
        \item For every pair of distinct points $x,y\in B(p,r)$,
        \[
        \frac{1}{d_Z(x',y')}\leq C_1\max\left\{\left(\frac{r}{d_Z(x,y)}\right)^\alpha,\left(\frac{r}{d_Z(x,y)}\right)^{1/\alpha}\right\}.
        \]
        
        \item For every $x\in B(p,r/2)$,
        \[
        B\!\left(x',\frac{1}{C_1}\right)\subset g\bigl(B(p,r)\bigr).
        \]
        
        \item If $F_\lambda:=Z\setminus B(p,\lambda r)$, then
        \[
        \operatorname{diam}\bigl(g(F_\lambda)\bigr) \leq \frac{\operatorname{diam}(Z)^2}{\delta} \eta_{\mathrm{qm}}\!\left(\frac{2}{\lambda}\right).
        \]
        In particular, $\operatorname{diam}(g(F_\lambda))\to0$ as $\lambda\to\infty$.
    \end{enumerate}
\end{lemma}

\begin{proof}
    By cocompactness of the action on $\operatorname{Pair}(Z)$ and Lemma~\ref{lem:cocompact_uniformly_separated}, there exists $\delta>0$ such that every pair of distinct points $u,v$ in $Z$ , there exists \(g \in G\) satisfying $d_Z(g(u),g(v)) \geq \delta$.

    Fix $p\in Z$ and $0<r<\operatorname{diam}(Z)$. Since $Z$ is connected, We may choose $q\in Z$ such that $d_Z(p,q)=r/2$. By the preceding separation property, there is $g\in G$ such that $d_Z(p',q')\geq\delta$. We prove the five assertions for this choice of $g$.

    \medskip
    \noindent\textbf{(1)} Let $x,y\in Z$ be distinct. Since $d_Z(p,q)=r/2$, at least one of $p,q$, say $a$, satisfies
    \[
    d_Z(x,a)\geq r/4.
    \]
    If $a\ne y$, power quasisymmetry gives
    \[
    \frac{d_Z(x',y')}{d_Z(x',a')} \leq C\max\left\{\left(\frac{d_Z(x,y)}{d_Z(x,a)}\right)^\alpha,\left(\frac{d_Z(x,y)}{d_Z(x,a)}\right)^{1/\alpha}\right\}.
    \]
    Using $d_Z(x',a')\leq\operatorname{diam}(Z)$ and $d_Z(x,a)\geq r/4$, we obtain the desired estimate after enlarging a constant depending only on $C$, $\alpha$, and $\operatorname{diam}(Z)$. If $a=y$, then the same conclusion follows directly from $d_Z(x',y')\leq\operatorname{diam}(Z)$ and $d_Z(x,y)\geq r/4$ after a further enlargement of the constant.

    \medskip
    \noindent\textbf{(2)} Since $g$ is a homeomorphism and $\operatorname{diam}(Z)>r$, there is a point $b\in Z\setminus\{x,y\}$ such that $d_Z(x',b')>r/2$. Indeed, the set of points at distance greater than $r/2$ from $x'$ is nonempty and open, and $Z$ has no isolated points. Applying power quasisymmetry to the triple $(x,b,y)$ gives
    \[
    \frac{d_Z(x',b')}{d_Z(x',y')} \leq C\max\left\{\left(\frac{d_Z(x,b)}{d_Z(x,y)}\right)^\alpha,\left(\frac{d_Z(x,b)}{d_Z(x,y)}\right)^{1/\alpha} \right\}.
    \]
    Since $d_Z(x,b)\leq\operatorname{diam}(Z)$ and $d_Z(x',b')>r/2$, the required estimate follows.

    \medskip
    \noindent\textbf{(3)} Let $x,y\in B(p,r)$ be distinct. Since $d_Z(p',q')\geq\delta$, at least one point $c\in\{p,q\}$ satisfies $d_Z(x',c')\geq\delta/2$. Also, $d_Z(x,c)\leq2r$. If $c\ne y$, power quasisymmetry yields
    \[
    \frac{\delta}{2d_Z(x',y')} \leq C\max\left\{\left(\frac{2r}{d_Z(x,y)}\right)^\alpha,\left(\frac{2r}{d_Z(x,y)}\right)^{1/\alpha}\right\}.
    \]
    This gives the asserted estimate after adjusting the constant. If $c=y$, then $d_Z(x',y')\geq\delta/2$ and $d_Z(x,y)\leq2r$, so the same estimate follows by enlarging the constant.

    \medskip
    \noindent\textbf{(4)} Let $x\in B(p,r/2)$ and $y\in Z\setminus B(p,r)$. Then $d_Z(x,y)\geq r/2$. Choose $c\in\{p,q\}$ such that $d_Z(x',c')\geq\delta/2$. Since $d_Z(x,c)\leq r$ and $c\in B(p,r)$, we have $c\ne y$. Power quasisymmetry gives
    \begin{align*}
        \frac{\delta/2}{d_Z(x',y')}
        &\leq\frac{d_Z(x',c')}{d_Z(x',y')}\\
        &\leq C\max\left\{\left(\frac{d_Z(x,c)}{d_Z(x,y)}\right)^\alpha,\left(\frac{d_Z(x,c)}{d_Z(x,y)}\right)^{1/\alpha}\right\}\\
        &\leq C\max\{2^\alpha,2^{1/\alpha}\}.
    \end{align*}
    Consequently,
    \[
    d_Z(x',y') \geq \rho:=\frac{\delta}{2C\max\{2^\alpha,2^{1/\alpha}\}}.
    \]
    Thus the ball $B(x',\rho)$ is disjoint from $g(Z\setminus B(p,r))$, and hence $B(x',\rho)\subset g(B(p,r))$. Choosing $C_1$ large enough that $1/C_1\leq\rho$ proves (4).

    \medskip
    \noindent\textbf{(5)} Let $y,z\in F_\lambda$ be distinct, and write
    \[
    s:=d_Z(y,p),\qquad t:=d_Z(z,p).
    \]
    Then $s,t\geq\lambda r$, while $d_Z(p,q)=r/2$. By the triangle inequality,
    \[
    d_Z(y,z)\leq s+t, \qquad d_Z(z,q)\geq t-r/2.
    \]
    It follows that
    \begin{align*}
        \frac{d_Z(y,z)d_Z(p,q)}{d_Z(y,p)d_Z(z,q)}
        &\leq\frac{(s+t)(r/2)}{s(t-r/2)}\\
        &=\frac r2\left(\frac1{t-r/2}
            +\frac{t}{s(t-r/2)}\right)\\
        &\leq\frac{1}{\lambda-1/2}\leq\frac{2}{\lambda},
    \end{align*}
    where the last inequality uses $\lambda\geq2$. Applying the uniform quasi-M\"obius inequality to the quadruple $(y,z,p,q)$ gives
    \[
    \frac{d_Z(y',z')d_Z(p',q')}{d_Z(y',p')d_Z(z',q')} \leq\eta_{\mathrm{qm}}\!\left(\frac{2}{\lambda}\right).
    \]
    Since $d_Z(p',q')\geq\delta$ and each distance in the denominator is at most $\operatorname{diam}(Z)$, we obtain
    \[
    d_Z(y',z')\leq\frac{\operatorname{diam}(Z)^2}{\delta}\eta_{\mathrm{qm}}\!\left(\frac{2}{\lambda}\right).
    \]
    Taking the supremum over all distinct $y,z\in F_\lambda$ proves (5). The right-hand side tends to zero as $\lambda\to\infty$.
\end{proof}

\section{Weak Tangents and Their One-Point Compactifications}
\label{sec:weak_tangents_one_point_compactification}

A \emph{pointed metric space} is a pair $(Z,p)$, where $(Z,d_Z)$ is a metric
space and $p\in Z$.

\begin{definition}[Pointed Gromov--Hausdorff convergence]
    A sequence $(Z_k,p_k)$ of pointed metric spaces is said to converge to a pointed metric space $(Z,p)$ if, for every $R>0$ and every $\varepsilon>0$, there exist $N\in\mathbb N$, a subset $M\subseteq B_Z(p,R)$, subsets $M_k\subseteq B_{Z_k}(p_k,R)$, and bijections $f_k:M_k\to M$ such that, for every $k\geq N$, the following conditions hold:
    \begin{enumerate}[$(i)$]
    \item $p\in M$, $p_k\in M_k$, and $f_k(p_k)=p$;
    \item $M$ is $\varepsilon$-dense in $B_Z(p,R)$, and $M_k$ is
    $\varepsilon$-dense in $B_{Z_k}(p_k,R)$;
    \item
    \[
    \bigl|d_{Z_k}(x,y)-d_Z(f_k(x),f_k(y))\bigr|<\varepsilon
    \qquad\text{for all }x,y\in M_k.
    \]
\end{enumerate}
\end{definition}

\begin{definition}[Weak tangent] \label{defn_weak_tangent}
    A complete pointed metric space $(S,q)$ is called a \emph{weak tangent} of a metric space $Z$ if there are numbers $\alpha_k>0$ with $\alpha_k\to\infty$ and points $p_k\in Z$ such that 
    \[
    (\alpha_k Z,p_k)\longrightarrow (S,q)
    \]
    in the pointed Gromov--Hausdorff sense. Here, for $\lambda>0$, we write $\lambda Z$ for the metric space $(Z,\lambda d_Z)$; thus, the underlying set is unchanged and its metric is multiplied by $\lambda$.
\end{definition}

The geometry of $Z$ is reflected in its weak tangents. Suppose $(Z,d_Z)$ is a compact metric space and $S\in\mathrm{WT}(Z)$ is a weak tangent of $Z$. If $Z$ is uniformly perfect, then $S$ is unbounded; if $Z$ is doubling, then $S$ is proper; and if $Z$ is Ahlfors $Q$-regular, then so is $S$ \cite{DS1997}.

\begin{lemma}\label{lem_weak_tangents_euclidean}
Every weak tangent of the standard round sphere $\mathbb S^n$ is pointed
isometric to $(\mathbb R^n,0)$. The same statement holds for
$\mathbb R^n$ equipped with its Euclidean metric.
\end{lemma}

\begin{proof}
For $\mathbb R^n$, translations followed by dilations identify
$(\alpha_k\mathbb R^n,p_k)$ isometrically with $(\mathbb R^n,0)$ for every
$k$. For the round sphere, rotations allow us to move each base point $p_k$
to a fixed point. In normal coordinates at that point, the rescaled
Riemannian metrics converge locally to the Euclidean metric. Thus every
pointed blow-up limit is $(\mathbb R^n,0)$.
\end{proof}

\paragraph{\textbf{One Point Compactification}}

Let $(S,d_S)$ be an unbounded locally compact metric space, and fix a base point
$q\in S$. Since $S$ is locally compact and noncompact, it has a one-point
compactification
\[
    \widehat S=S\cup\{\infty\}.
\]
To define a metric on $\widehat S$ associated with the pointed space
$(S,q)$, let $h_q:\widehat S\to[0,\infty)$ be given by
\[
 h_q(x)=
 \begin{cases}
 \dfrac{1}{1+d_S(x,q)},&x\in S,\\[6pt]
 0,&x=\infty.
 \end{cases}
\]
For $x,y\in S$, set
\[
\rho_q(x,y)=h_q(x)h_q(y)d_S(x,y),
\]
and, for $x\in S$, set
\[
    \rho_q(x,\infty)=\rho_q(\infty,x)=h_q(x),
    \qquad \rho_q(\infty,\infty)=0.
\]
The function $\rho_q$ is symmetric and separates points, but it need not
satisfy the triangle inequality. We therefore replace it by its associated
chain metric. For $x,y\in\widehat S$, define
\[
    \widehat d_q(x,y) :=\inf\left\{ \sum_{i=0}^{m-1}\rho_q(x_i,x_{i+1}): m\in\mathbb N,\ x_0=x,\ x_m=y,\ x_i\in\widehat S \right\}.
\]

\begin{lemma}[{\cite[Lemma~2.2]{BK2002QM}}] \label{lem_metric_one_point_compactification}
    Let $(S,d_S)$ be an unbounded locally compact metric space, let $q\in S$, and let $\widehat S=S\cup\{\infty\}$ be its one-point compactification. Define $\rho_q$ and $\widehat d_q$ as above. Then:
    \begin{enumerate}[(1)]
        \item For all $x,y\in\widehat S$,
        \[
        \frac14\,\rho_q(x,y)\leq \widehat d_q(x,y)\leq\rho_q(x,y).
        \]
        
        \item The function $\widehat d_q$ is a metric on $\widehat S$, and its induced topology agrees with the one-point compactification topology. In particular, $(\widehat S,\widehat d_q)$ is compact. 
        
        \item The identity map
        \[
        \operatorname{id}_S:(S,d_S)\longrightarrow (S,\widehat d_q|_S)
        \]
        is a $16t$-quasi-M\"obius homeomorphism.
    \end{enumerate}
\end{lemma}

\begin{proposition}[{\cite[Proposition~4.4]{PS2024}}] \label{prop_weak_tangent_compactification}
    Suppose that $Z$ is a compact, uniformly perfect, doubling metric space, and that $G\curvearrowright Z$ is a uniformly quasisymmetric action for which the induced action $G\curvearrowright\operatorname{Pair}(Z)$ is cocompact. If $(S,q)\in\mathrm{WT}(Z)$, then there exists a quasisymmetric homeomorphism
    \[
    h:(\widehat S,\widehat d_q)\longrightarrow Z.
    \]
    Moreover, the restriction
    \[
    h|_S:S\longrightarrow Z\setminus\{h(\infty)\}
    \]
    is a quasi-M\"obius homeomorphism.
\end{proposition}

\begin{lemma} \label{lem:WTS_implies_QM_spaces}
    Suppose that $X$ and $Y$ are compact, uniformly perfect, doubling metric spaces. Suppose also that $G\curvearrowright X$ and $H\curvearrowright Y$ are uniformly quasi-M\"obius actions for which the induced actions on $\operatorname{Pair}(X)$ and $\operatorname{pair}(Y)$, respectively, are cocompact. If there are weak tangents $(S,p)\in\mathrm{WT}(X)$ and $(T,q)\in\mathrm{WT}(Y)$ and a quasi-M\"obius homeomorphism $f:S\to T$, then there exists a quasi-M\"obius homeomorphism $g:X\to Y$. In particular, since $X$ and $Y$ are bounded, $g$ is also quasisymmetric.
\end{lemma}

\begin{proof}
    Let $(S,p)\in\mathrm{WT}(X)$ and $(T,q)\in\mathrm{WT}(Y)$, and let $f:S\to T$ be a quasi-M\"obius homeomorphism. Since $f$ is a homeomorphism between locally compact spaces, it extends to a homeomorphism between their one-point compactifications by setting
    \[
    \widehat f(x)=
    \begin{cases}
       f(x),&x\in S,\\
       \infty,&x=\infty.
    \end{cases}
    \]
    By Lemma~\ref{lem_metric_one_point_compactification}, the identity maps from $(S,d_S)$ to $(S,\widehat d_p|_S)$ and from $(T,d_T)$ to $(T,\widehat d_q|_T)$ are quasi-M\"obius. Thus, on the complements of $\infty$, we can write 
    \[
    \widehat f|_S =\operatorname{id}_T\circ f\circ\operatorname{id}_S^{-1}, 
    \]
    where the identity maps are understood between the corresponding metrics. Since inverses and compositions of quasi-M\"obius homeomorphisms are quasi-M\"obius, $\widehat f|_S$ is quasi-M\"obius between the punctured compactifications. The quasi-M\"obius inequality extends to quadruples containing $\infty$ by approximating $\infty$ with a sequence in $S$ and using continuity of the compactification metrics. Thus $\widehat f$ is quasi-M\"obius on the whole compactifications. This is the same compactification argument used in the proof of \cite[Lemma~5.8]{BK2002QM}.

    By Proposition~\ref{prop_weak_tangent_compactification}, there are quasi-M\"obius homeomorphisms \[
    h_X:(\widehat S,\widehat d_p)\longrightarrow X, \qquad h_Y:(\widehat T,\widehat d_q)\longrightarrow Y. 
    \]
    Define
    \[
    g=h_Y\circ\widehat f\circ h_X^{-1}:X\longrightarrow Y. 
    \] 
    The inverse and the composition of quasi-M\"obius homeomorphisms are quasi-M\"obius. Hence $g$ is a quasi-M\"obius homeomorphism. Since $X$ is bounded, it follows that $g$ is quasisymmetric as well. 
\end{proof}

\section{Stable Values, Bounded Multiplicity, and Regular Maps}
\label{sec:stable_values_regular_maps}

In this section, we recall several notions and results from
\cite{BK2002QM} concerning continuous maps between metric spaces. They
will be used in the sequel.

\begin{definition}[Stable value]
    Let $X$ and $Y$ be metric spaces, and let $f:X\to Y$ be a continuous map. A point $y\in Y$ is called a \emph{stable value} of $f$ if there exists $\varepsilon>0$ such that $y\in\operatorname{Im}(g)$ for every continuous map $g:X\to Y$ satisfying
    \[
    \operatorname{dist}(f,g)<\varepsilon,
    \]
    where
    \[
    \operatorname{dist}(f,g)
    :=\sup_{x\in X}d_Y\bigl(f(x),g(x)\bigr)
    \]
    is the uniform distance between $f$ and $g$.
\end{definition}

A continuous map is called \emph{light} if each of its fibres is totally disconnected. The following proposition provides a sufficient condition for a continuous map to have a stable value.

\begin{proposition}[{\cite[Proposition~3.2]{BK2002QM}}] \label{prop:light_maps_stable_values}
    Let $X$ be a compact metric space with topological dimension $\dim_{\mathrm{top}}X\geq n$, where $n\in\mathbb N$, and let $f:X\to\mathbb R^n$ be a light continuous map. Then $f$ has a stable value.
\end{proposition}

We next recall the notion of bounded multiplicity and a structural
result for continuous maps satisfying this condition.

\begin{definition}[Bounded multiplicity]
    A map $f:X\to Y$ between two spaces is said to have \emph{bounded multiplicity} if there exists $N\in\mathbb N$ such that 
    \[
    \#f^{-1}(y)\leq N \qquad\text{for every }y\in Y.
    \]
\end{definition}

\begin{theorem}[{\cite[Theorem~3.4]{BK2002QM}}] \label{thm:bounded_multiplicity_covering}
    Let $X$ be a compact metric space such that every nonempty open subset of $X$ has topological dimension at least $n$. Suppose that $f:X\to\mathbb R^n$ is a continuous map of bounded multiplicity. Then there exists an open set $V\subseteq\mathbb R^n$ such that 
    \[
    V\subseteq\operatorname{Im}(f) \qquad\text{and}\qquad \overline V=\operatorname{Im}(f),
    \]
    where the closure is taken in $\mathbb R^n$. Moreover, setting $U:=f^{-1}(V)$, the set $U$ is dense in $X$, and the restriction 
    \[
    f|_U:U\longrightarrow V 
    \]
    is a covering map.
\end{theorem}

We finally recall a notion that will be useful when studying maps
between weak tangents.

\begin{definition}[Regular map]
    Let $X$ be a metric space and let $Y$ be a doubling metric space. A map $f:X\to Y$ is called \emph{regular} if it is Lipschitz and there exists a constant $N\in\mathbb N$ such that, for every open ball $B_Y(y,r)\subseteq Y$, the preimage $f^{-1}(B_Y(y,r))$ can be covered by at most $N$ open balls of radius $r$ in $X$.
\end{definition}

\begin{remark}
    Every regular map has bounded multiplicity. Indeed, suppose that there are $N+1$ distinct points $x_1,\ldots,x_{N+1}\in f^{-1}(y)$ for some $y\in Y$. Set
    \[
    \delta:=\min\{d_X(x_i,x_j):i\neq j\}>0.
    \]
    By regularity, the set $f^{-1}(B_Y(y,\delta/2))$ can be covered by at most $N$ open balls of radius $\delta/2$ in $X$. Each such ball contains at most one of the points $x_1,\ldots,x_{N+1}$, since its diameter is at most $\delta$ and the balls are open. This is a contradiction. Hence $\#f^{-1}(y)\leq N$ for every $y\in Y$.
\end{remark}

The following proposition shows that a Lipschitz map with an image of positive $Q$-dimensional Hausdorff measure induces a regular map between suitable weak tangents.

\begin{proposition}[{\cite[Proposition~4.1]{BK2002QM}}] \label{prop:regular_map_between_weak_tangents}
    Let $f:X\to Y$ be a Lipschitz map, where $X$ is a compact Ahlfors $Q$-regular metric space for some $Q>1$, and $Y$ is a complete doubling metric space. If $\mathcal H^Q\bigl(f(X)\bigr)>0$, then there exist weak tangents $S\in\mathrm{WT}(X)$ and $T\in\mathrm{WT}(Y)$, together with a regular map $g:S\to T$. 
\end{proposition}

\section{Proof of Theorem~\ref{thm_main_1}}
\label{sec:proof_main_1}

Throughout this section, let $n\in\mathbb N$, and let $Z$ be a compact, Ahlfors $n$-regular metric space with $\dim_{\mathrm{top}}(Z)=n$. Suppose that $G\curvearrowright Z$ is a uniformly quasi-M\"obius action for which the induced action $G\curvearrowright\operatorname{Pair}(Z)$ is cocompact.

The group of M\"obius transformations of $\mathbb S^n$ acts uniformly quasi-M\"obiusly, and its induced action on $\operatorname{Pair}(\mathbb S^n)$ is cocompact. Hence, by Lemma~\ref{lem:WTS_implies_QM_spaces}, it suffices to show that some weak tangent of $Z$ is quasisymmetrically equivalent to $\mathbb R^n$.

\begin{theorem} \label{thm:tangent_quasisymmetric_to_euclidean}
    Let $n\in\mathbb N$, and let $Z$ be a compact, Ahlfors $n$-regular metric space with $\dim_{\mathrm{top}}Z=n$. Suppose that $G\curvearrowright Z$ is a uniformly quasi-M\"obius action and that the induced action $G\curvearrowright\operatorname{Pair}(Z)$ is cocompact. Then there exists a weak tangent $S$ of $Z$ that is quasisymmetrically equivalent to $\mathbb R^n$.
\end{theorem}

Before proving Theorem~\ref{thm:tangent_quasisymmetric_to_euclidean}, we
establish the following lemmas.

\begin{lemma} \label{lem_top_dim_of_open_subsets}
    Suppose that $Z$ is a compact, uniformly perfect and doubling metric space, and that $G\curvearrowright Z$ is a uniformly quasi-M\"obius action for which the induced action $G\curvearrowright\operatorname{Pair}(Z)$ is cocompact. If $\dim_{\mathrm{top}}(Z)=n\in\mathbb N$, then $\dim_{\mathrm{top}}(U)=n$ for every nonempty open subset $U$ of $Z$ or of any weak tangent of $Z$.
\end{lemma} 

\begin{proof}
    We first consider a nonempty open subset $U\subseteq Z$. Choose a nonempty open ball $B=B(p,r)$ such that $\overline B\subseteq U$. We claim that there is a sequence $\{g_k\}$ in $G$ such that \begin{equation}\label{eq:shrinking_complement_of_gB} 
        \operatorname{diam}\bigl(Z\setminus g_k(\overline B)\bigr) \longrightarrow 0.
    \end{equation}
    Indeed, since $Z$ is uniformly perfect, we can choose a sequence $\lambda_k\to0$ and points
    \[
    x_k^1=p,\qquad x_k^2\in B(p,\lambda_k r)\setminus B\left(p,\frac{\lambda_k r}{C}\right), 
    \] 
    where $C>1$ is independent of $k$. Thus $d_Z(x_k^1,x_k^2)> (\lambda_k/C)r$. By the cocompactness of the action on $\operatorname{Pair}(Z)$, there is a constant $\delta>0$ and, for every $k$, an element $g_k\in G$ such that
    \[
    d_Z\bigl(g_k(x_k^1),g_k(x_k^2)\bigr)\geq\delta.
    \]
    Since $Z$ is compact, a uniformly quasi-M\"obius family of self-homeomorphisms of $Z$ is uniformly quasisymmetric. We may therefore apply part~(2) of Lemma~\ref{lem:tech_2} to the balls $B_k=B(p,r)$ and the points $x_k^1,x_k^2$. It follows that
    \[
    \operatorname{diam}\bigl(Z\setminus g_k(B)\bigr)\longrightarrow0.
    \]
    Since $g_k(B)\subseteq g_k(\overline B)$, this implies \eqref{eq:shrinking_complement_of_gB}. It follows that the set
    \[
    E:=Z\setminus\bigcup_{k\in\mathbb N}g_k(\overline B)
    \]
    contains at most one point. Each $g_k(\overline B)$ is closed in $Z$, and $E$ is either empty or a singleton, hence is also closed. By the countable sum theorem for topological dimension (see \cite[Theorem~II.1]{Nagata1983}), and the invariance of topological dimension under homeomorphisms, we obtain 
    \[
    \begin{aligned}
        \dim_{\mathrm{top}}(Z) 
        &\leq \sup\left\{ \dim_{\mathrm{top}}\bigl(g_k(\overline B)\bigr), \dim_{\mathrm{top}}(E)\right\}\\
        &\leq \max\left\{\dim_{\mathrm{top}}(\overline B),0\right\}.
    \end{aligned}
    \]
    Since $\overline B\subseteq U$ and topological dimension is monotone for subspaces of metric spaces, this gives
    \[
    n=\dim_{\mathrm{top}}(Z) \leq\dim_{\mathrm{top}}(U) \leq\dim_{\mathrm{top}}(Z)=n.
    \]
    Thus $\dim_{\mathrm{top}}(U)=n$. 
    
    Now let $U$ be a nonempty open subset of a weak tangent $S$ of $Z$. By Proposition~\ref{prop_weak_tangent_compactification}, the one-point compactification of $S$ is homeomorphic to $Z$. In particular, the restriction of this homeomorphism maps $S$ homeomorphically onto the complement of one point in $Z$. Since $S$ is open in its one-point compactification, the image of $U$ is a nonempty open subset of $Z$. The first part of the proof therefore gives $\dim_{\mathrm{top}}(U)=n$. 
\end{proof}

\begin{lemma} \label{lem:lipschitz_map_positive_measure}
    Suppose that $Z$ is as in Theorem~\ref{thm:tangent_quasisymmetric_to_euclidean}. Then there exists a Lipschitz map $h:Z\to\mathbb S^n$ such that 
    \[
    \mathcal H^n\bigl(h(Z)\bigr)>0.
    \]
\end{lemma}

\begin{proof}
    Since $\dim_{\mathrm{top}}Z=n$, by \cite[Theorem~III.1]{Nagata1983} there exists a continuous map $g:Z\to\mathbb S^n$ with a stable value $y\in\mathbb S^n$.
    Let $g_1,\ldots,g_{n+1}$ be the coordinate functions of $g$, where $\mathbb S^n$ is viewed as the unit sphere in $\mathbb R^{n+1}$. Since $Z$ is a compact metric space, Lipschitz functions are uniformly dense in $C(Z)$. Hence, for each $i$, we can choose a Lipschitz function $f_i:Z\to\mathbb R$ that uniformly approximates $g_i$ as closely as desired. Thus the map 
    \[ 
    f=(f_1,\ldots,f_{n+1}):Z\longrightarrow\mathbb R^{n+1} 
    \]
    is Lipschitz and can be chosen arbitrarily close to $g$. 
    
    Choose the approximation sufficiently close that $\|f(z)\|\geq 1/2$ for every $z\in Z$, and define
    \[
    \pi:\mathbb R^{n+1}\setminus\{0\}\longrightarrow\mathbb S^n, \qquad \pi(v)=\frac{v}{\|v\|},
    \]
    followed by $h:=\pi\circ f$. The map $h$ is Lipschitz. Moreover, by making the approximation sufficiently close, we can ensure that $h$ is uniformly close to $g$.

    Let $\varepsilon_0>0$ be such that every continuous map $u:Z\to\mathbb S^n$ with $\operatorname{dist}(u,g)<\varepsilon_0$ has $y$ in its image. Choose $f$ so that $\operatorname{dist}(h,g)<\varepsilon_0/2$. It follows that $y$ is also a stable value of $h$: if $\operatorname{dist}(u,h)<\varepsilon_0/2$, then $\operatorname{dist}(u,g)<\varepsilon_0$, and hence $y\in u(Z)$.

    We next observe that a stable value of a map into $\mathbb S^n$ lies in the interior of its image. Indeed, if $y$ were not an interior point of $h(Z)$, there would be a sequence $y_j\in\mathbb S^n\setminus h(Z)$ converging to $y$. Choose rotations $R_j$ of $\mathbb S^n$ such that $R_j(y_j)=y$ and $R_j\to\operatorname{id}_{\mathbb S^n}$ uniformly. Indeed, the angle of this rotation is $\theta_j=d_{\mathbb S^n}(y_j,y)$, which tends to zero, consequently, the rotations approach the identity uniformly: $\sup_{x\in\mathbb S^n}d_{\mathbb S^n}(R_j(x),x)\to 0$. Then $R_j\circ h$ converges uniformly to $h$, while $y\notin R_j(h(Z))$, contradicting the stability of $y$ for $h$. Consequently, $h(Z)$ contains a nonempty open subset of $\mathbb S^n$, and therefore
    \[
    \mathcal H^n\bigl(h(Z)\bigr)>0.
    \]
\end{proof}

\begin{proof}[Proof of Theorem~\ref{thm:tangent_quasisymmetric_to_euclidean}]
    By Lemma~\ref{lem:lipschitz_map_positive_measure} and Proposition~\ref{prop:regular_map_between_weak_tangents}, there exist a weak tangent $S$ of $Z$, a weak tangent $T$ of $\mathbb S^n$, and a regular map
    \[
    \phi:S\longrightarrow T.
    \]
    Every weak tangent of the standard sphere $\mathbb S^n$ is isometric to $\mathbb R^n$, so we identify $T$ with $\mathbb R^n$. Since $\phi$ is regular, it has bounded multiplicity. Moreover, by Lemma~\ref{lem_top_dim_of_open_subsets}, every nonempty open subset of $S$ has topological dimension $n$.

    Since $Z$ is compact and Ahlfors $n$-regular, it is uniformly perfect and doubling, and so $S$ is proper and unbounded. Choose a bounded nonempty open subset $W\subseteq S$, and set $X:=\overline W$. Properness implies that $X$ is compact. Every nonempty relatively open subset $A$ of $X$ meets $W$ in a nonempty open subset of $S$. Hence
    \[
    \dim_{\mathrm{top}}(A) \geq \dim_{\mathrm{top}}(A\cap W)=n.
    \]
    We may therefore apply Theorem~\ref{thm:bounded_multiplicity_covering} to the restriction $\phi|_X:X\to\mathbb R^n$. We obtain an open subset $V\subseteq\phi(X)$ such that $\overline V=\phi(X)$, and, on setting $U_0:=(\phi|_X)^{-1}(V)$, the set $U_0$ is dense in $X$ and
    \[
    \phi|_{U_0}:U_0\longrightarrow V
    \]
    is a covering map.

    Since $U_0$ is dense in $X$ and $W$ is a nonempty relatively open subset of $X$, choose $x_0\in U_0\cap W$. The covering map is locally a homeomorphism. Thus we may choose an open ball $B\subseteq V$ centered at $\phi(x_0)$ and an open subset $U\subseteq U_0\cap W$ containing $x_0$ such that the restriction
    \[
    \psi:=\phi|_U: U\longrightarrow B
    \]
    is a homeomorphism. 
    By \cite[Lemma~4.2]{BK2002QM}, the map $\psi$ is bi-Lipschitz. 

    Now fix $x\in U$ and put $y:=\psi(x)$. Applying \cite[Lemma~4.3]{BK2002QM} to the bi-Lipschitz map $\psi$ shows that $S$ has a weak tangent that is bi-Lipschitz equivalent to a weak tangent of $\mathbb R^n$. Since every weak tangent of $\mathbb R^n$ is isometric to $\mathbb R^n$, there is a weak tangent $S'$ of $S$ that is bi-Lipschitz equivalent to $\mathbb R^n$.

    Finally, weak tangents are transitive: since $S$ is a weak tangent of $Z$ and $S'$ is a weak tangent of $S$, the space $S'$ is a weak tangent of $Z$. Thus $Z$ has a weak tangent that is bi-Lipschitz, and hence quasisymmetrically, equivalent to $\mathbb R^n$. This completes the proof.
\end{proof}

\section{De-snowflaking a Metric and Proof of Theorem~\ref{thm_main_2}}
\label{sec:desnowflaking_metric}

In this section, we construct a metric \(d_{\mathrm{new}}\) that is uniformly comparable to \(d^{1/\epsilon}\), under suitable geometric and quantitative assumptions. Thus, the original metric is comparable to the \(\epsilon\)-snowflake of \(d_{\mathrm{new}}\).

\subsection{Snowflake maps and regularity}

Let $(X,d_X)$ and $(Y,d_Y)$ be metric spaces. A map $f:X\to Y$ is called a \emph{snowflake map with exponent $\alpha>0$} if there is a constant $C\geq 1$ such that 
\begin{equation}\label{eq:snowflake-map}
    C^{-1}d_X(x,y)^\alpha \leq d_Y(f(x),f(y)) \leq C d_X(x,y)^\alpha \qquad (x,y\in X).
\end{equation}
In particular, $f$ is injective and is a homeomorphism onto its image. If $\alpha=1$, it is bi-Lipschitz onto its image. If $0<\alpha<1$, then small nonzero distances are enlarged relative to the original metric, whereas if $\alpha>1$, they are diminished. 
Every snowflake map is quasisymmetric, with control function \(\eta(t)=C^2t^\alpha\).

Suppose that $(X,d,\mu)$ is Ahlfors $Q$-regular, and let $0<\alpha<1$. The function $d_\alpha(x,y)=d(x,y)^\alpha$ is a metric, since $(s+t)^\alpha\leq s^\alpha+t^\alpha$ for $s,t\geq0$. Moreover, $B_{d_\alpha}(x,r)=B_d(x,r^{1/\alpha})$, and hence
\[
\mu(B_{d_\alpha}(x,r)) =\mu(B_d(x,r^{1/\alpha})) \asymp r^{Q/\alpha}.
\]
Consequently, \((X,d^\alpha,\mu)\) is Ahlfors
\(Q/\alpha\)-regular. 

For \(\alpha>1\), the function \(d^\alpha\) need not satisfy the triangle inequality. For example, on \(\mathbb R\) with its usual metric,
\[
d_\alpha(0,2)=2^\alpha>2 =d_\alpha(0,1)+d_\alpha(1,2).
\]
More generally, for $\alpha>0$, whenever \(d^\alpha\) is a metric, the Hausdorff dimensions satisfy \begin{equation}\label{eq:hausdorff-dimension-snowflake}
    \dim_H(X,d^\alpha) =\frac{1}{\alpha}\dim_H(X,d),
\end{equation}
and $(Z,d^\alpha)$ is a $Q/\alpha$-Ahlfors regular.

\subsection{Proof of Theorem~\ref{thm_main_2}}

As is Theorem~\ref{thm_main_2}, let \(n \in \mathbb{N}\) with $n\geq 2$, \(0 < \epsilon \le 1\), and let \(Z = (Z,d)\) be a compact metric space, homeomorphic to \(S^n\), and Ahlfors regular of dimension \(Q = n/\varepsilon\). Suppose that \(\Gamma \curvearrowright Z\) is a uniformly power quasisymmetric action with a control function $\eta(t)= C \max\{t^\alpha,t^{1/\alpha}\}$ for some $C \geq 1$ and $\alpha>0$, such that the induced action on $\operatorname{Pair}(Z)$ is cocompact. Assume, moreover, that \(Z\) satisfies the following discrete length property: 
\begin{equation}\label{eq:discrete-length-property}
    \text{each } \tau\text{-path between two points } x,y \text{ has length} \gtrsim (d(x,y)/{\tau})^{1/\epsilon}
\end{equation}

We claim that there exist a metric \(d_{\mathrm{new}}\) on \(Z\) and a constant \(C_*\geq1\) such that
\begin{equation}\label{eq:prop41}
    C_*^{-1}d(x,y)^{1/\epsilon} \leq d_{\mathrm{new}}(x,y) \leq C_*d(x,y)^{1/\epsilon}
\end{equation}
for all \(x,y\in Z\).

\paragraph{\textbf{Definition of the new metric}}

Fix \(a>2\). For each \(k\in\mathbb N\), let \(S_k\) be a maximal \(a^{-k}\)-separated subset of \(Z\). Then \(S_k\) is an \(a^{-k}\)-net. Define the collection of \(k\)-balls by
\[
\mathcal B_k=\{B(s,2a^{-k}):s\in S_k\}.
\]
A finite sequence \(B_1,\ldots,B_\ell\) of \(k\)-balls is called a \emph{\(k\)-ball chain} if \(B_i\cap B_{i+1}\neq\varnothing\) for \(1\leq i<\ell\). It connects \(x,y\in Z\) if \(x\in B_1\) and \(y\in B_\ell\). Its length is the number \(\ell\) of balls, counted with multiplicity. Since \(Z\) is compact and connected, any two points can be joined by a finite \(k\)-ball chain. For \(x,y\in Z\), set
\[
d_k(x,y) = \bigl(\text{minimum length of a \(k\)-ball chain joining \(x\) to \(y\)}\bigr) a^{-k/\epsilon}.
\]
The function \(d_k\) is symmetric and satisfies the triangle inequality, but \(d_k(x,x)=a^{-k/\epsilon}\), so it is not a metric. Define
\[
d_{\mathrm{new}}(x,y)=\limsup_{k\to\infty}d_k(x,y).
\]
We will prove uniform estimates
\begin{equation}\label{eq:dk-two-sided}
    d(x,y)^{1/\epsilon}\lesssim d_k(x,y) \lesssim d(x,y)^{1/\epsilon}
\end{equation}
for all distinct \(x,y\in Z\) and all sufficiently large \(k\). The implicit constants will be independent of \(x,y\), and \(k\). 

Once these estimates are established, taking the limsup in \eqref{eq:dk-two-sided} as \(k\to\infty\) yields
\[
d(x,y)^{1/\epsilon} \lesssim d_{\mathrm{new}}(x,y) \lesssim d(x,y)^{1/\epsilon}.
\]
Finally, \(d_{\mathrm{new}}(x,x)=0\), and it is symmetric. The triangle inequality for each \(d_k\) implies
\[
d_{\mathrm{new}}(x,z) \leq d_{\mathrm{new}}(x,y)+d_{\mathrm{new}}(y,z),
\]
since all the limsups are finite by the established upper bound. The lower bound shows that \(d_{\mathrm{new}}(x,y)>0\) whenever \(x\neq y\). Thus \(d_{\mathrm{new}}\) is a metric and satisfies \eqref{eq:prop41} as desired.

\medskip

\paragraph{\textbf{Lower bound of (\ref{eq:dk-two-sided}) estimation}} 
Let \(x,y\in Z\) be distinct, and consider a \(k\)-ball chain \(B_1,\ldots,B_\ell\) joining them. For each \(1\leq i<\ell\), choose \(z_i\in B_i\cap B_{i+1}\), and put \(z_0=x\), \(z_\ell=y\). Since every \(k\)-ball has diameter at most \(4a^{-k}\), we have
\[
d(z_{i-1},z_i)\leq4a^{-k} \qquad (1\leq i\leq\ell).
\]
Applying \eqref{eq:discrete-length-property} with \(\tau=4a^{-k}\), we obtain
\[
\ell\geq C_D^{-1} \left(\frac{d(x,y)}{4a^{-k}}\right)^{1/\epsilon}.
\]
Multiplying by \(a^{-k/\epsilon}\) and taking a shortest chain gives
\[
d_k(x,y)\geq C_D^{-1}4^{-1/\epsilon}d(x,y)^{1/\epsilon}.
\]
This proves the lower bound in \eqref{eq:dk-two-sided}.

\medskip

\paragraph{\textbf{Upper bound of (\ref{eq:dk-two-sided}) estimation}}

We establish the following chain-connecting estimate, which will play a key role in deriving the upper bound in \eqref{eq:dk-two-sided}.

\begin{proposition} \label{prop:ball-chain}
    Let $Z$ is a metric space as in Theorem~\ref{thm_main_2}. Let \(x,y \in Z\) and \(m \in \mathbb{N}\) with $a^{- m} < d(x,y) \le a^{-(m-1)}$. Then for each \(k \ge m\), there is a \(k\)-ball chain connecting the balls \(B_x=B(x,a^{- m})\) and \(B_y=B(y,a^{- m})\) of length at most \(C_2 a^{(k-m)/\epsilon}\). Here, \(C_2\) depends only on the data associated to \(Z\).
\end{proposition}

We now use Proposition~\ref{prop:ball-chain} to establish the desired upper bound in \eqref{eq:dk-two-sided}. The proof of the proposition is postponed to a later subsection. 

\begin{proof}[Proof of the upper bound]
    Let \(x,y\in Z\) be distinct, and choose \(h\in\mathbb N\) such that 
    \[
    a^{-h}<d(x,y)\leq a^{-(h-1)}.
    \] 
    Fix \(k >> h\) sufficiently large. We denote
    \[
    B_m^z = B(z,a^{-m}), \qquad \text{for \(z\in Z\) and \(m\in\mathbb{N}\)}.
    \]
    
    At the first step, by applying Proposition~\ref{prop:ball-chain}, we obtain a \(k\)-ball chain joining \(B_h^x\) and \(B_h^y\) of length at most \(C_2a^{(k-h)/\epsilon}\). This gives us points
    \[
    x = x_{1,0},\, x_{1,1},\, x_{1,2},\, x_{1,3} = y
    \]
    where \(x_{1,0},x_{1,1}\in B_h^x\) and \(x_{1,2},x_{1,3}\in B_h^y\), and \(x_{1,1}\) is connected to \(x_{1,2}\) by a \(k\)-ball chain of length at most \(C_2a^{(k-h)/\epsilon}\).
    
    At the second step, observe that \(d(x_{1,0},x_{1,1})\le a^{-h}\) and so we can apply the Proposition~\ref{prop:ball-chain} again to obtain a \(k\)-ball chain of length at most \(C_2a^{(k-h-1)/\epsilon}\) connecting \(B_{h+1}^{x_{1,0}}\) and \(B_{h+1}^{x_{1,1}}\). This gives us points
    \[
    x = x_{2,0},\, x_{2,1},\, x_{2,2},\, x_{2,3} = x_{1,1}
    \]
    where \(x_{2,0},x_{2,1}\in B_{h+1}^{x_{1,0}}\) and \(x_{2,2},x_{2,3}\in B_{h+1}^{x_{1,1}}\), and \(x_{2,1}\) is connected to \(x_{2,2}\) by a \(k\)-ball chain of length at most \(C_2a^{(k-h-1)/\epsilon}\). A similar process will be apply to the pair of points \(x_{1,2}\) and \(x_{1,3}\).

    For \(j=2,\ldots,k-h+1\), proceed inductively. Suppose that at step \(j-1\) we have \(2^j\) ordered points
    \[
    x = x_{j-1,0},x_{j-1,1},\ldots,x_{j-1,2^j-1}=y,
    \]
    such that
    \[
    d(x_{j-1,i},x_{j-1,i+1}) \leq a^{-(h+j-2)}
    \]
    for every even \(i\) with \(0\leq i\leq 2^j-2\), and there are already constructed \(k\)-ball chains joining \(x_{j-1,i+1}\) to \(x_{j-1,i+2}\) for every even \(i\) with \(0\leq i\leq 2^j-4\). For each such even \(i\), use the preceding consequence with \(m=h+j-1\) to choose 
    \[
    x_{j,2i+1}\in B(x_{j-1,i},a^{-(h+j-1)}), \qquad x_{j,2i+2}\in B(x_{j-1,i+1},a^{-(h+j-1)}) 
    \]
    so that these two new points are joined by a \(k\)-ball chain of length at most \(C_2 a^{(k-h-j+1)/\epsilon}\). Retain the old points in the order
    \[
    x_{j,2i}=x_{j-1,i}\quad(i\text{ even}),\qquad x_{j,2i+1}=x_{j-1,i}\quad(i\text{ odd}), 
    \]
    and insert the new points in the remaining positions. This produces \(2^{j+1}\) ordered points. 

    At stage \(j\), exactly \(2^{j-1}\) new chains are introduced, each of length at most \(C_2 a^{(k-h-j+1)/\epsilon}\), while the short gaps between consecutive points in even positions have size at most \(a^{-(h+j-1)}\). The previously constructed chains join the intervening odd-even pairs. 

    After \(k-h+1\) stages, we obtain \(2^{k-h+2}\) ordered points
    \[
    x=x_{k-h+1,0},\ldots, x_{k-h+1,2^{k-h+2}-1}=y,
    \]
    such that the distance between consecutive points in each even position is at most \(a^{-k}\). Each such pair is contained in a single \(k\)-ball: choose a net point within distance \(a^{-k}\) of the first point, and use the fact that the second is at distance at most \(a^{-k}\) from the first. Concatenating these single \(k\)-balls with all the previously constructed chains gives a \(k\)-ball chain from \(x\) to \(y\) of length at most
    \[
    2^{k-h+1} +\sum_{j=1}^{k-h+1} 2^{j-1} C_2 a^{(k-h-j+1)/\epsilon} \leq \left(2+\frac{C_2}{1-2a^{-1/\epsilon}}\right) a^{(k-h)/\epsilon} \lesssim a^{(k-h)/\epsilon}.  
    \] 
    Here we used \(a^{1/\epsilon}>2\), and $a$ is fixed. Then the implicit constant depends only on \(C_2\) and \(\epsilon\). Therefore 
    \[
    d_k(x,y) \lesssim a^{(k-h)/\epsilon}a^{-k/\epsilon} = a^{-h/\epsilon} \leq d(x,y)^{1/\epsilon}.
    \]
    This completes the proof of upper bound of \eqref{eq:dk-two-sided}.
\end{proof}

\subsection{Proof of Proposition~\ref{prop:ball-chain}}

Now we prove the chain-connecting estimate Proposition~\ref{prop:ball-chain}. We use the following standard separation lemma for the sphere. 

We first observe the packing estimation that for each \(k\), the collection of all \(k\)-balls that meet \(B(p,r)\) has an upper bound. 

\begin{lemma}\label{lem:size_of_sep_sets}
    If \(Z\) is Ahlfors Ahlfors \(Q\)-regular, then there exists \(C_P\) depending on $Q$ such that, whenever \(0<t\leq r\), every \(t\)-separated subset \(S\subset B(p,r)\) satisfies
    \begin{equation}\label{eq:local-packing}
        \#S\leq C_P\left(\frac r t\right)^Q.
    \end{equation}
\end{lemma}

\begin{proof}
    Let \(S=\{x_1,\ldots,x_N\}\subset B(z,r)\) be a \(t\)-separated set, meaning $d(x_i,x_j)\geq t$ for $i\neq j$. Then the balls \(B(x_i,t/2)\) are pairwise disjoint. By Ahlfors \(Q\)-regularity, we have $\mu(B(x_i,t/2))\gtrsim (t/2)^Q$ for each $i=1,\dots,N$. Since \(S \subset B(z,r)\), we have
    \[
    \bigcup_{i=1}^N B(x_i,t/2)
    \subset B(z,r+t/2).
    \]
    Consequently, 
    \[
    N\cdot(t/2)^Q \lesssim \sum_{i=1}^N\mu(B(x_i,t/2)) \leq \mu(B(z,r+t/2)) \lesssim (3r/2)^Q,
    \]
    Therefore,
    \[
    N\lesssim \left(\frac rt\right)^Q.
    \]
\end{proof}

\begin{lemma}[Lemma~4.5]\label{lem:separated-balls}
    There is an absolute constant \(C_3>0\) with the following property. Let \(B_0,B_1\) be metric balls of radius \(\delta>0\) in \(\mathbb S^n\) such that \(\operatorname{dist}(B_0,B_1)\geq\delta\). Suppose that \(E\subset\mathbb S^n\) satisfies 
    \[
    \operatorname{diam}(E)<\delta, \qquad \operatorname{dist}(B_i,E)\geq\delta \quad (i=0,1). 
    \]
    Then there is a topological \(n\)-cube \(\widetilde S\subset\mathbb S^n\setminus E\) with a pair of opposite codimension-one faces contained in \(B_0\) and \(B_1\), respectively, such that any two points lying in opposite codimension-one faces of \(\widetilde S\) are at spherical distance at least \(C_3\delta^3\).
\end{lemma}

\begin{proposition}[Chain-product inequality, Proposition~4.4]
\label{prop:chain-product}
    Let \(\mathcal U\) be a finite open cover of a topological \(n\)-cube, and let \(F_i,G_i\), \(1\leq i\leq n\), be its \(n\) pairs of opposite codimension-one faces. Let \(\ell_i\) be the minimum number of members of \(\mathcal U\) in a chain joining \(F_i\) to \(G_i\), where consecutive members of a chain are required to intersect. Then 
    \[
    \#\mathcal U\geq\ell_1\ell_2\cdots\ell_n.
    \]
\end{proposition}

\begin{proof}[Proof of Proposition~\ref{prop:ball-chain}]
    Let \(x,y \in Z\) and \(m \in \mathbb{N}\) with $a^{- m} < d(x,y) \le a^{-(m-1)}$. We claim that there exists $C_2>0$ for every \(k\geq m\) there is a \(k\)-ball chain connecting \(B(x,a^{-m})\) to \(B(y,a^{-m})\) whose length is at most
    \begin{equation}\label{eq:ball-chain-bound}
        C_2a^{(k-m)/\epsilon}.
    \end{equation}

    Let $\beta = \max\{\alpha,1/\alpha\}$, and choose \(C'\geq 2\) sufficiently large so that it larger that the constant $C_1$ considered in Lemma~\ref{lem:power_quasisymmetric_estimates}. Fix \(x,y,m\) as in the claim and set 
    \[
    p=x,\qquad r=2C' a^{-(m-1)}, \qquad C''= C'^{-1}(2C'a)^{\beta}.
    \]
    Then \(x,y\in B(p,r/2C')\), and \(B(x,a^{-m}),B(y,a^{-m})\subset B(p,r)\). 

    Fix a homeomorphism \(f:Z\to\mathbb S^n\). By uniform continuity of \(f^{-1}\), choose \(\delta>0\) sufficiently small and $\lambda \geq 2$ large enough so that 
    \begin{align} 
        d_Z(u,v)\geq \frac{1}{C''} 
        &\ \Longrightarrow\ d_{\mathbb S^n}(f(u),f(v))\geq3\delta, \label{eq:sphere-ball-inclusion-1}\\ 
        d_Z(u,v)<\omega(1/\lambda) 
        &\ \Longrightarrow\ d_{\mathbb S^n}(f(u),f(v))<\frac{C_3\delta^3}{2}. \label{eq:sphere-ball-inclusion-2} 
    \end{align}

    By applying Lemma~\ref{lem:power_quasisymmetric_estimates}, we have \(g \in G\) \(p,r\) corresponding to the above fixed $p$, $r$ and $\lambda$ such the the followings are hold. Write \(x'=g(x)\) and \(y'=g(y)\).

    Since \(d(x,y)>a^{-m}\) and \(d(x,y)\leq a^{-(m-1)}\), property~$(iii)$ of Lemma~\ref{lem:power_quasisymmetric_estimates} gives
    \[
    \frac{1}{d(x'y')} \leq C_1 \max \left\{\left( \frac{r}{d(x,y)}\right)^\alpha, \left(\frac{r}{d(x,y)}\right)^{1/\alpha}\right\} \leq C_1 \left(\frac{r}{d(x,y)}\right)^\beta.
    \]
    Implies,
    \begin{equation}\label{eq:7.10.3}
        d(x',y') 
        \geq \frac{1}{C_1} \left(\frac{a^{-m}}{2C' a^{-(m-1)}}\right)^\beta 
        = \frac{1}{2^\beta C'^{(\beta+1)} a^\beta}\\
        \geq \frac{1}{C''}
    \end{equation}

    We next show that
    \begin{equation}\label{eq:small-image-balls}
        B(x',1/C'')\subset g(B(x,a^{-m})), \qquad B(y',1/C'')\subset g(B(y,a^{-m})). 
    \end{equation} 
    For the first inclusion, let \(w\in B(x',1/C'')\). Since \(B(x',1/C'') \subset g(B(p,r))\), the point \(z=g^{-1}(w)\) belongs to \(B(p,r)\). Again, by applying property~$(iii)$ of Lemma~\ref{lem:power_quasisymmetric_estimates} to \(p,z\), and using \(d(p,z)<r\), yields
    \[
    d(x,z) \leq r\left(\frac{1}{C_1}\,d(g(p),g(z))\right)^{1/\beta} =r\left(\frac{1}{C_1} \,d(x',w) \right)^{1/\beta} < 2C'a^{-(m-1)}(\frac{1}{C_1C''})^{1/\beta} \leq a^{-m}. 
    \]
    Thus \(z\in B(x,a^{-m})\), proving the first inclusion. The second inclusion follows in the same way. 

    Now property~$(iv)$ of Lemma~\ref{lem:power_quasisymmetric_estimates} gives  \begin{equation}\label{eq:4.7}
        B(x',1/C_1) \cup B(y',1/C_1) \subset g(B(p,r)),
    \end{equation} 
    and if we denote $F_\lambda = Z \setminus B(p,\lambda r)$, then as \(g\) is a homeomorphism, this implies that
    \begin{equation}\label{eq:7.10.5}
        \operatorname{dist}\bigl(\{x',y'\}, g(F_\lambda)\bigr) \ge \frac{1}{C'} \geq \frac{1}{C''}.
    \end{equation}
    
    Moreover, property~$(v)$ of Lemma~\ref{lem:power_quasisymmetric_estimates} gives    
    \begin{equation}
        \operatorname{diam}(g(F_\lambda)) \leq \omega(1/\lambda), \qquad \text{and} \qquad \omega(1/\lambda) \to 0 \text{ as } \lambda\to\infty.
    \end{equation}

    \vspace{1cm}    
    Let us now use the homeomorphism \(f : Z \to \mathbb{S}^n\). By applying \eqref{eq:7.10.3}, \eqref{eq:7.10.5} and our choice of \(\delta\) from \eqref{eq:sphere-ball-inclusion-1}, we have
    \[
    d_{S^n}\bigl(f(x'),f(y')\bigr) \ge 3\delta \quad \text{and} \quad \operatorname{dist}\bigl(\{f(x'),f(y')\},fg(F_\lambda)\bigr) \ge 3\delta 
    \]
    We denote, 
    \[
    B_0 = B_{\mathbb{S}^n}(f(x'),\delta), \qquad B_1 = B_{\mathbb{S}^n}(f(y'),\delta), \qquad E= fg(F_\lambda).
    \]
    Then $B_0$ and $B_1$ are metric balls in \(\mathbb{S}^n\) with $\operatorname{dist}(B_0,B_1) \ge \delta$ and of distance at least \(\delta\) from the set \(E \subset \mathbb{S}^n\). Thus $B_0$, $B_1$ and $E$ satisfy the hypotheses of Lemma~\ref{lem:separated-balls}. 
    Furthermore, \eqref{eq:sphere-ball-inclusion-1} and \eqref{eq:small-image-balls} imply 
    \[
    B_0\subset f(g(B(x,a^{-m}))), \qquad B_1\subset f(g(B(y,a^{-m}))). 
    \]

    Apply Lemma~\ref{lem:separated-balls} to obtain a topological \(n\)-cube \(\widetilde S\subset\mathbb S^n\setminus E\), with opposite faces \(\widetilde X_1\subset B_0\) and \(\widetilde Y_1\subset B_1\), such that points in any pair of opposite codimension-one faces are at spherical distance at least \(C_3 \delta^3\). Set 
    \[
    S=(f\circ g)^{-1}(\widetilde S). 
    \]
    Since \(\widetilde S\cap E=\varnothing\), we have \(S\subset B(p,\lambda r)\). The set \(S\) is a topological \(n\)-cube. The corresponding opposite faces 
    \[
    X_1=(f\circ g)^{-1}(\widetilde X_1), \qquad Y_1=(f\circ g)^{-1}(\widetilde Y_1) 
    \]
    satisfy 
    \[
    X_1\subset B(x,a^{-m}), \qquad Y_1\subset B(y,a^{-m}). 
    \]

    Fix \(k\geq m\), and let \(\mathcal U\) be the collection of all \(k\)-balls that meet \(S\). This is a finite open cover of \(S\). Since \(S\subset B(p,\lambda r)\), the packing estimate Lemma~\ref{lem:size_of_sep_sets} gives
    \begin{equation}\label{eq:cover-cardinality}
        \#\mathcal U \lesssim \left(\frac{\lambda r}{a^{-k}} \right)^Q\lesssim a^{Q(k-m)}.
    \end{equation}
    For \(1\leq i\leq n\), let \(\ell_i\) be the minimum number of members of \(\mathcal U\) in a chain joining the \(i\)-th pair \(X_i,Y_i\) of opposite codimension-one faces of \(S\). 
    We claim that, for each \(i=2,\ldots,n\),
    \begin{equation}\label{eq:other-face-chain-lower}
        \ell_i\gtrsim a^{(k-m)/\epsilon}.
    \end{equation}
    Indeed, consider any chain of \(\ell_i\) \(k\)-balls joining \(X_i\) to \(Y_i\). Choosing points in successive intersections gives a discrete \(4a^{-k}\)-path in \(Z\) from some \(z\in X_i\) to some \(w\in Y_i\), with at most \(\ell_i\) steps. Since \(f(g(z))\) and \(f(g(w))\) lie in opposite faces of \(\widetilde S\), we have 
    \[
    d_{\mathbb S^n}(f(g(z)),f(g(w)))\geq C_3\delta^3.
    \]
    By the choice of \(lamb
    \) in \eqref{eq:sphere-ball-inclusion-2}, this implies
    \[
    d(g(z),g(w))\geq\ \omega(1/\lambda).
    \]
    Property~$(i)$ of Lemma~\ref{lem:power_quasisymmetric_estimates} now gives
    \[
    \omega(1/\lambda) \leq C_1 \max \left\{ \left( \frac{d(z,w)}{r} \right)^\alpha, \left( \frac{d(z,w)}{r} \right)^{1/\alpha}  \right\}
    = C_1 \left( \frac{d(z,w)}{r} \right)^{\beta'},
    \]
    where $\beta' \in \{\alpha,1/\alpha\}$. Since $\alpha\lambda$ is fixed, so is $\beta',\omega(1/\lambda)$. Then 
    \[
    d(z,w) \geq \left( \frac{1}{C_1} \omega(1/\lambda) \right)^{1/\beta} r \gtrsim r \gtrsim a^{-m}.
    \]
    Therefore, the discrete length property \eqref{eq:discrete-length-property} yields 
    \[
    \ell_i\gtrsim \left(\frac{a^{-m}}{4a^{-k}}\right)^{1/\epsilon} =a^{(k-m)/\epsilon}, 
    \]
    which proves \eqref{eq:other-face-chain-lower}. 

    Applying Proposition~\ref{prop:chain-product} to \(\mathcal U\), and using \(Q=n/\epsilon\), we obtain 
    \[
    a^{Q(k-m)} \gtrsim \#\mathcal U \geq \ell_1\ell_2\cdots\ell_n  \gtrsim \ell_1\,a^{(n-1)(k-m)/\epsilon}. 
    \]
    It follows that
    \[
    \ell_1\lesssim a^{\left(Q-(n-1)/\epsilon\right)(k-m)} =a^{(k-m)/\epsilon}.
    \]
    A chain realizing \(\ell_1\) joins \(X_1\) to \(Y_1\), and hence joins \(B(x,a^{-m})\) to \(B(y,a^{-m})\). This proves \eqref{eq:ball-chain-bound}.
\end{proof}

\section{Proof of Theorem~\ref{thm:generalized-rel-Cannon-conjec}}
\label{sec:proof_main_3}

As in Theorem~\ref{thm:generalized-rel-Cannon-conjec}, let \((G,\mathcal P)\) be a relatively hyperbolic group pair whose Bowditch boundary \(Z=\partial(G,\mathcal P)\) is homeomorphic to \(\mathbb S^n\), with \(n\geq 2\). Assume that \(Z\) admits an Ahlfors \(Q\)-regular metric \(d\), where \(Q=n/\epsilon\) for some \(0<\epsilon\leq 1\). Furthermore, suppose there exists a constant \(c>0\) such that, for every pair of points \(x,y\in Z\) and every discrete \(\tau\)-path joining them, its combinatorial length \(\ell\) satisfies
\[
\ell\geq c\left(\frac{d(x,y)}{\tau}\right)^{1/\epsilon}.
\]

Since \((G,\mathcal P)\) is relatively hyperbolic and its Bowditch boundary is homeomorphic to \(\mathbb S^n\), the action \(G\curvearrowright Z\) is a non-elementary geometrically finite convergence action, by uniformly power quasisymmetry homeomorphisms. Moreover, the induced action on the space of distinct pairs \(\operatorname{Pair}(Z)\) is cocompact. 

    By the discrete length hypothesis and Theorem~\ref{thm_main_2}, there exists a metric \(d_{\mathrm{new}}\) on \(Z\) such that
    \[
    d_{\mathrm{new}}(x,y)\asymp d(x,y)^{1/\epsilon}.
    \]
    Moreover, \((Z,d_{\mathrm{new}})\) is Ahlfors \(n\)-regular and belongs to the same quasisymmetric gauge as \((Z,d)\), see Corollary~\ref{cor:thm-2}. Since the original action is uniformly power quasisymmetric, the conjugated action on \((Z,d_{\mathrm{new}})\) is uniformly quasisymmetric. Cocompactness on pairs is also preserved, since \(d\) and \(d_{\mathrm{new}}\) induce the same topology.
    Consequently, Theorem~\ref{thm_main_1} gives a quasisymmetric homeomorphism
    \[
    h:(Z,d_{\mathrm{new}})\longrightarrow \mathbb S^n.
    \]
    Conjugating the action by \(h\), we obtain an action
    \[
    \rho:G\longrightarrow \operatorname{Homeo}(\mathbb S^n),
    \qquad
    \rho(g)=h\circ g\circ h^{-1}.
    \]
    Since \(h\) is quasisymmetric and the action on \((Z,d_{\mathrm{new}})\) is uniformly quasisymmetric, the group \(\rho(G)\) is uniformly quasisymmetric on the standard sphere. By the equivalence between quasisymmetric and quasiconformal mappings on \(\mathbb S^n\), the action \(\rho(G)\) is uniformly quasiconformal. 

    We now apply Tukia's theory of uniformly quasiconformal groups. By Tukia's invariant conformal structure theorem (Theorem~F in \cite{Tukia1986}), there exists a bounded, measurable, \(\rho(G)\)-invariant conformal structure 
    \[
    \mu\colon \mathbb{S}^n \longrightarrow \mathcal{S}(n) = \mathrm{SL}(n,\mathbb{R})/\mathrm{SO}(n).
    \]
    We claim that \(\mu\) is approximately continuous at some conical limit point of the \(\rho(G)\)-action. Since \(\mathcal{S}(n)\), equipped with its Riemannian metric, is a separable metric space, the Lebesgue differentiation theorem implies that \(\mu\) is approximately continuous almost everywhere on \(\mathbb{S}^n\); see \cite{Heinonen2001,EV1992}. Let \(E_1 \subseteq \mathbb{S}^n\) denote the set of points at which \(\mu\) is approximately continuous. Thus, \(E_1\) has full spherical Lebesgue measure.
    

    On the other hand, since \(G\) is countable and \(\mathcal{P}\) is finite, the set of parabolic points of the Bowditch boundary \(Z\) is countable, being a finite union of \(G\)-orbits. By geometric finiteness, every point of \(Z\) is either conical or bounded parabolic. Consequently, the set of conical limit points has full measure with respect to any measure that assigns zero measure to countable sets. In particular, under the homeomorphism \(h\colon Z \to \mathbb{S}^n\), the set of parabolic points remains countable. Hence, the set \(E_2 \subseteq \mathbb{S}^n\) of conical limit points of the \(\rho(G)\)-action has full spherical Lebesgue measure.   

     It follows that \(E_1 \cap E_2\) has full spherical Lebesgue measure and, in particular, is nonempty. Thus, there exists a conical limit point at which \(\mu\) is approximately continuous, proving the claim. 

    We can now apply Tukia's rigidity theorem (Theorem~G in \cite{Tukia1986}). Since \(\rho(G)\) is uniformly quasiconformal, preserves the bounded measurable conformal structure \(\mu\), and \(\mu\) is approximately continuous at the conical limit point \(\xi\), there exists a quasiconformal homeomorphism $f:\mathbb S^n\to\mathbb S^n$ such that
    \[
    f\rho(G)f^{-1}\subseteq \operatorname{\text{M\"ob}}(\mathbb S^n),
    \]
    where \(\operatorname{\text{M\"ob}}(\mathbb S^n)\) denotes the group of M\"obius transformations of the standard sphere. It follows that the composition \(f\circ h:Z\to\mathbb S^n\) conjugates the original action of \(G\) on \(Z\) to an action by M\"obius transformations. Let
\[
\Phi=f\circ h:Z\longrightarrow\mathbb S^n, \qquad \text{and} \qquad \Gamma=f\rho(G)f^{-1}\leq\operatorname{\text{M\"ob}}(\mathbb S^n).
\]

\medskip
We now show that the action of \(G\) on \(\mathbb H^{n+1}\) is properly discontinuous by isometries. Since \(G\curvearrowright Z\) is a convergence action, the induced action on the space of distinct triples \(\operatorname{Tri}(Z)\) is properly discontinuous. The homeomorphism \(\Phi:Z\to\mathbb S^n\) conjugates this action to the action of \(\Gamma\) on \(\operatorname{Tri}(\mathbb S^n)\), which is therefore also properly discontinuous. On the other hand, the M\"obius group \(\operatorname{\text{M\"ob}}(\mathbb S^n)\) acts properly on \(\operatorname{Tri}(\mathbb S^n)\). It follows that \(\Gamma\) is a discrete subgroup of \(\operatorname{\text{M\"ob}}(\mathbb S^n)\). Under the standard identification $\operatorname{\text{M\"ob}}(\mathbb S^n)\cong\operatorname{Isom}^{+}(\mathbb H^{n+1})$, the action of \(\Gamma\) extends to a properly discontinuous isometric action on \(\mathbb H^{n+1}\). Consequently, \(G\) admits a properly discontinuous isometric action on \(\mathbb H^{n+1}\).

\medskip
Finally, we show that every peripheral subgroup is virtually \(\mathbb Z^n\). Let \(P\in\mathcal P\) be a peripheral subgroup, and let \(\xi\in Z\) be its associated bounded parabolic point. Set \(\xi'=\Phi(\xi)\in\mathbb S^n\). Then \(\Gamma_P=\operatorname{Stab}_{\Gamma}(\xi')\) is the image of \(P\) under the conjugated action.

After conjugating by a M\"obius transformation, we may assume that \(\xi'=\infty\) in the upper half-space model of \(\mathbb H^{n+1}\). Every element of \(\Gamma_P\) fixes \(\infty\). Since \(\xi\) is a bounded parabolic point, its stabilizer contains no loxodromic elements. Thus, every element of \(\Gamma_P\) is either parabolic or elliptic. Consequently, each \(g\in\Gamma_P\) has the form
\[
g(x,t)=(A_gx+b_g,t),
\]
where \(A_g\in O(n)\) and \(b_g\in\mathbb R^n\). In particular, \(\Gamma_P\) preserves the horosphere \(H=\{(x,1):x\in\mathbb R^n\}\) and acts on it by Euclidean isometries.

Since \(\Gamma\) is discrete in \(\operatorname{Isom}(\mathbb H^{n+1})\), the induced action of \(\Gamma_P\) on \(H\cong\mathbb R^n\) is properly discontinuous. Moreover, since \(\xi\) is a bounded parabolic point, the action \(P\curvearrowright Z\setminus\{\xi\}\) is cocompact. By the conjugacy \(\Phi\), the induced action of \(\Gamma_P\) on \(\mathbb S^n\setminus\{\xi'\}\cong\mathbb R^n\) is cocompact. Under the identification of \(\mathbb S^n\setminus\{\xi'\}\) with the horosphere \(H\), this is precisely the induced Euclidean action. Hence \(\Gamma_P\) acts properly discontinuously and cocompactly by Euclidean isometries on \(\mathbb R^n\). By the first Bieberbach theorem, \(\Gamma_P\) contains a finite-index subgroup isomorphic to \(\mathbb Z^n\). Therefore, every peripheral subgroup \(P\in\mathcal P\) is virtually \(\mathbb Z^n\). This completes the proof.

\bigskip
\paragraph{\textbf{Declaration on the Use of Generative AI}}
The authors used the free version of ChatGPT (GPT-6, OpenAI) solely for language editing, including grammatical corrections and improvements to the clarity and readability of the manuscript. The authors take full responsibility for the content of the manuscript. 

\bibliographystyle{amsalpha}
\bibliography{Bibliography} 


\end{document}